\documentclass[11pt]{amsart}
\usepackage{graphicx}
\usepackage{amssymb}

\usepackage{amssymb,amsmath,amsthm}
\usepackage[all]{xy}
\usepackage{url}
  \usepackage{color}
  \usepackage{xy}
 \xyoption{curve}
 \xyoption{color}
 \xyoption{line}
\xyoption{arc}
\usepackage{graphicx}
\usepackage[
           colorlinks=true,
            urlcolor=blue,       
            filecolor=blue,     
            linkcolor=blue,       
            citecolor=blue,         
            pdftitle= {Title},
            pdfauthor={Author},
            pdfsubject={Subject},
            pdfkeywords={Key}
            pagebackref,
            pdfpagemode=None,
            bookmarksopen=true]
            {hyperref}
\usepackage{hyperref}
\usepackage{color}
\usepackage[author=,english,status=draft]{fixme}
\usepackage{ifdraft}

\ifdraft{
\usepackage{fancyhdr}
\pagestyle{fancy}
\fancyhead{}
\fancyfoot{}
\fancyhead[CO,CE]{Draft --- \today}
\fancyhead[RO, LE] {\thepage}
}

\date{}
\dedicatory{To Miles Reid}

 \newtheorem{thm}{Theorem}[section]

 \theoremstyle{definition}
 \newtheorem{example}[thm]{Example}
 
 \newtheorem{remark}[thm]{Remark}
\numberwithin{equation}{section}

\newcommand{\bbC}{{\mathbb{C}}}

\newcommand{\bbG}{{\mathbb{G}}}

\newcommand{\bbP}{{\mathbb{P}}}

\newcommand{\bbQ}{{\mathbb{Q}}}
\newcommand{\bbT}{{\mathbb{T}}}
\newcommand{\bbZ}{{\mathbb{Z}}}
\newcommand{\bbF}{{\mathbb{F}}}

\newcommand{\bfF}{{\mathbf{F}}}

\newcommand{\Foc}{\operatorname{Foc}}

\newcommand{\SL}{\operatorname{SL}}

\newcommand{\PGL}{\operatorname{PGL}}

\newcommand{\rank}{\operatorname{rank}}
\newcommand{\Aut}{\operatorname{Aut}}

\newcommand{\CR}{\operatorname{CR}}
\newcommand{\Kum}{\operatorname{Kum}}

\newcommand{\half}{\frac{1}{2}}
\newcommand{\da}{\dasharrow}
\newcommand{\St}{\operatorname{St}}
\newcommand{\st}{\operatorname{st}}
\newcommand{\Bit}{\operatorname{Bit}}
\newcommand{\Gr}{\operatorname{Gr}}

\newcommand{\bsm}{\left(\begin{smallmatrix}}
\newcommand{\esm}{\end{smallmatrix}\right)}

\newcommand{\la}{\langle}
\newcommand{\ra}{\rangle}

\newcommand{\calC}{\mathcal{C}}
\newcommand{\calE}{\mathcal{E}}

\newcommand{\calH}{\mathcal{H}}

\newcommand{\calL}{\mathcal{L}}
\newcommand{\calM}{\mathcal{M}}
\newcommand{\calN}{\mathcal{N}}
\newcommand{\calO}{\mathcal{O}}
\newcommand{\calP}{\mathcal{P}}
\newcommand{\calQ}{\mathcal{Q}}

\newcommand{\calS}{\mathcal{S}}

\newcommand{\calD}{\mathcal{D}}

\newcommand{\frakr}{\mathfrak{r}}

\newcommand{\frakS}{\mathfrak{S}}

\newcommand{\sfD}{\mathsf{D}}

\newcommand{\sfS}{\mathsf{S}}

\newcommand{\sfV}{\mathsf{V}}
\newcommand{\dP}{\operatorname{dP}}

\newcommand{\Pic}{\mathrm{Pic}}
\newcommand{\Fano}{\mathrm{Fano}}

\newcommand{\Rey}{\operatorname{Rey}}
\newcommand{\rey}{\operatorname{Rey}}
\newcommand{\Seg}{\operatorname{Seg}}

\newcommand{\beq}{\begin{equation}}

\newcommand{\eeq}{\end{equation}}

\title[15-nodal quartic surfaces.I]{15-nodal quartic surfaces.I: quintic del Pezzo surfaces and congruences of lines in $\bbP^3$}
\author{Igor V. Dolgachev}
\address{Department of Mathematics,
University of Michigan, Ann Arbor, 48109}
\email{idolga@umich.edu}

\begin{document}
\begin{abstract}
We explain a classical construction of a del Pezzo surface of degree 
$d = 4$ or $5$ as a smooth order $2$ congruence of lines in $\bbP^3$ whose focal surface is a quartic surface $X_{20-d}$ with $20-d$  ordinary double points. We also show that $X_{15}$ can be realized  as  a hyperplane section of the Castelnuovo-Richmond-Igusa quartic hypersurface in $\bbP^4$. This leads to the proof  of rationality of the moduli space of $15$-nodal quartic surfaces. We discuss some other birational models of $X_{15}$: quartic symmetroids,  $5$-nodal quartic surfaces, $10$-nodal sextic surfaces in $\bbP^4$ and nonsingular surfaces of degree 10 in $\bbP^6$. Finally we study some birational involutions of a 15-nodal quartic surface which, as it is shown in  Part 2 of the paper jointly with I. Shimada, belong to a finite set of  generators of the  group of birational automorphisms of a general $15$ nodal quartic  surface.
\end{abstract}

\maketitle
\tableofcontents
\section{Introduction}

The geometry of 16-nodal Kummer quartic surfaces is widely discussed in classical as well as in modern literature. This paper addresses less known but nevertheless very rich geometry of 15-nodal quartic surfaces. What unites the two classes of nodal quartic surfaces is the fact that surfaces from both classes are realized in 6 different ways as focal surfaces of smooth congruences of lines in $\bbP^3$, congruences of bidegree $(2,2)$ in the former case and bidegree $(2,3)$ in the latter case. The line geometry realization of the quartic surfaces brings an unexpected dividend: it defines a certain involution on a smooth K3 model of such surface with quotient a del Pezzo surface of degree 4 in the Kummer case and degree 5 in the other case. Other interesting involutions come from a realization of the surface as a quartic symmetroid, the discriminant surface of a web of quadrics in $\bbP^3$. The quotient by these involutions is a Coble surface with the bi-anticanonical divisor consisting of six disjoint smooth rational curves in the Kummer case and five such curves in the other case. 

Of course, a Kummer surface comes as a specialization of a 15-nodal quartic surface by smoothing one of its 16 ordinary nodes. Explicitly, it can be seen by realizing any 15-nodal quartic surface as a hyperplane section of the notorious Castelnuovo-Richmond-Igusa quartic hypersurface in $\bbP^4$ with 15 double lines, the tangent hyperplanes give a Kummer surface. The fact that such  realization is always possible (but not generically possible as it is well known) which we prove here  (also proved independently by A. Avilov in \cite{Avilov}) is used to prove  the rationality of the moduli space of 15-nodal quartic surfaces. The relationship between $15$-nodal quartic surfaces $X_{15}$ and the beautiful geometry of the Castelnuovo-Richmond-Igusa quartic hypersurface and its dual Segre cubic hypersurface goes back to H.F. Baker \cite[Chapter V]{Baker}. We employ this relationship to describe some other interesting birational models of $X_{15}$ such as a $5$-nodal quartic surface with $20$ lines, a sextic model with $10$ nodes, a smooth model as a quadric section of the Fano threefold of genus 5 and index 2. 

The  first half of the paper is devoted to summarizing some classical results about  congruences of lines in $\bbP^3$ and especially congruences of order 2 from \cite{Jessop} and \cite{Sturm}. Following the classical approach we do not use the assumption of smoothness of the congruence assuming only isolated non-normal singularities leaving a possibility to extend some of our results to quartic surfaces with less than $15$ number of nodes.  

The  group of birational automorphisms a general 16-nodal Kummer quartic surface was described by S. Kondo \cite{Kondo}. In Part 2 \cite{Part2} of the paper (jointly with Ichiro Shimada) we  find a finite set of generators that contains some of the involutions discussed in the paper.

We work over $\bbC$, however most of the results which we discuss here are valid  when the ground field is an algebraically closed field of characteristic $p\ne 2,3$.

The paper originates from an attempt to answer some questions of S. Mukai related to $15$-nodal quartic surfaces. I am thankful to him for many discussions, questions and insights on this and on other closely related questions about geometry of Enriques and Coble surfaces. I am also grateful to A. Verra and I. Shimada for  helpful consultations.

\section{Generalities on congruences of lines}\label{S:generalities}
\textbf{2.1} Following classical terminology, a congruence of lines is a complete irreducible $2$-dimensional family of lines in $\bbP^3$. It is parameterized by a surface $S$ in  the Grassmannian variety $\bbG = \Gr_1(\bbP^3)$ of lines in $\bbP^3$. The lines $\ell_s$ corresponding to points  $s\in S$ are called \emph{rays} to distinguish them from any line in $\bbP^3$. 

The cohomology class $[S]$ of $S$ in $H^4(\bbG,\bbZ) \cong \bbZ^2$ is equal to $m[\sigma_x]+n[\sigma_\pi]$, where $\sigma_x$ is a plane in $\bbG$ of rays containing a point $x\in \bbP^3$ and $\sigma_\pi$ is a plane in $\bbG$ parametrizing lines in a plane $\pi\subset \bbP^3$. The number $m$ (resp. $n$) is called the \emph{order} (resp. \emph{class}) of $S$ and the pair $(m,n)$ is referred to as the bidegree of $S$. The sum $m+n$ is the degree of $S$ in the Pl\"ucker space $\bbP^5$. The natural  duality between lines in $\bbP^3$ and  lines in the dual $\bbP^3$ defines an isomorphism between the Grassmannians of lines, the image $S^*$ of $S$ under this isomorphism is called the \emph{dual congruence}. Its bidegree is $(n,m)$. For brevity, we will assume that $m,n\ne 1$ (a  modern treatment of other cases can be found in \cite{Ran}). 

 The universal family of lines $Z_S = \{(x,s)\in \bbP^3\times S:x\in \ell_s\}$ comes with two projections: 
$$\xymatrix{\bbP^3&Z_S\ar[l]_p\ar[r]^q&S.}$$ 
The restriction of $p$ to a fiber $q^{-1}(s)$ is an isomorphism onto the ray $\ell_s$. This allows us to identify the fibers of $q$ with the rays of the congruence. A confusing classical terminology defines a \emph{singular point} of $S$ as a point $x\in \bbP^3$ such that the fiber $p^{-1}(x)$ consists of infinitely many rays. A one-dimensional irreducible component of the set of singular points is called a \emph{fundamental curve} of $S$. 

 We will assume that $S$ has no fundamental curves and it is smooth or has only isolated singular points such that  the normalization $S'$ of $S$ is smooth. The rays corresponding to singular points of $S$ are called \emph{multi-rays}. This kind of singularities arise from projections of smooth surfaces $S'$ in $\bbP^6$ from a point that is contained in a finitely many secant lines to $S'$. If $S$ is not smooth and $\phi:S'\to S$ is the normalization map,  we replace $Z_S$ with the base change $Z_{S'}$, and get the following diagram 
$$\xymatrix{&Z_{S'}\ar[dl]_{p'}\ar[r]^{q'}\ar[d]&S'\ar[d]^\phi\\
\bbP^3&Z_S\ar[l]_p\ar[r]^q&S}.$$  
 Under these assumptions, $Z_{S'}$  is smooth and  the map $p'$ is a generically finite map of degree $m$ of smooth $3$-dimensional varieties. The ramification divisor of $p'$
 $$R_{S'}:= \textrm{Ram}(p')\subset Z_{S'}$$
  is a closed subscheme of $Z_{S'}$ equal to the support of $\Omega_{Z_{S'}/\bbP^3}^1$. The usual Hurwitz type formula gives
 \beq\label{canformula}
\omega_{Z_{S'}}  \cong p'{}^*\omega_{\bbP^3}\otimes \calO_{Z_{S'}}(R_{S'}) \cong p'{}^*\calO_{\bbP^3}(-4)\otimes \calO_{Z_{S'}}(R_{S'}).
\eeq 
The projection $q':Z_{S'}\to S'$ is the projectivization of the pre-image of the universal rank $2$ tautological vector bundle   over $\bbG$ under the map 
$\phi':S'\overset{\phi}{\to} S\hookrightarrow \bbG$. Its first Chern class is equal to $h:= c_1(\calO_{S'}(1))$, where $\calO_{S'}(1)$ is the pre-image of $\calO_{\bbG}(1)$ under the map $\phi'$. The usual formula for the canonical class of a projective bundle gives another formula for the canonical sheaf of $Z_{S'}$
\beq\label{omegaZ}
\omega_{Z_{S'}} \cong q'{}^*(\omega_{S'}\otimes \calO_{S'}(-1))\otimes p'{}^*(\calO_{\bbP^3}(-2)).
\eeq 
Comparing the two formulas, we find 
\beq\label{ramdiv}
[R_{S'}] = p'{}^*(2H)+q'{}^*(h+K_{S'}),
\eeq
where $H = c_1(\calO_{\bbP^3}(1))$. The image of $R_{S'}$ under the projection to $\bbP^3$ is the branch divisor of $p'$. It is called the \emph{focal surface} of $S$ and will be denoted by $\Foc(S)$.  Intersecting $R_{S'}$ with a general fiber of $q'$, we deduce  from  formula \eqref{ramdiv}   that 
\beq\label{degree2}
\deg q'|_{R_{S'}} = 2.
\eeq
If the fiber $q'{}^{-1}(s')$ is contained in $R_{S'}$, its image $\ell_s$ in $\bbP^3$ lies in $\Foc(S)$. Otherwise $q'{}^{-1}(s')$ contains two critical points in $R_{S'}$ and the ray $\ell_s$ intersects $\Foc(S)$ at the images of these points with multiplicity $\ge 2$.  It follows from \eqref{degree2} that  the focal surface $\Foc(S)$ consists of at most two irreducible components. We will assume that it is irreducible. We will also assume that $\Foc(S)$ contains only finitely many rays of $S$, so that the surface $R_{S'}$ is irreducible too.
 Hence
\begin{itemize}
\item[(1)]\it{ $S$ is an irreducible component of the surface of bitangents of $\Foc(S)$}.
\end{itemize}

For any curve $C$ on $S$, the image of $q^{-1}(C)$ under the map $p:Z_S\to \bbP^3$ is a  \emph{ruled surface}. Its degree is equal to the degree of $C$ in the Pl\"ucker embedding. For any  line $l$ in $\bbP^3$, the set of rays intersecting $l$ is a (special) hyperplane section $R(l) = \bbT_l(\bbG)\cap S$ of $S$. Following the classical notation we denote by $(l)$ the corresponding ruled surface in $\bbP^3$.

 Let $\calP_l\cong \bbP^1$ be the pencil of planes containing a general line $l$ in $\bbP^3$. Consider a map
$$f:R(l) \to \calP_l\times l \cong \bbP^1\times \bbP^1$$
that assigns  to a point $s\in R(l),$ the pair $(\la l,\ell_s\ra,\ell_s\cap l)$ that consists of the plane $\la l,\ell_s\ra$ generated by the line $l$ and the ray $\ell_s$ and the intersection point $\ell_s\cap l$.  The composition of $f$ with the first projection (resp. the second projection) is a map of degree $n$ (resp. $m$). For a general line $l$, the curve $R(l)$ is smooth and its genus is called the \emph{sectional genus} of $S$. The image of the map $f$ is  a curve of bidegree $(m,n)$ on $\bbP^1\times \bbP^1$ and arithmetic genus $(m-1)(n-1)$.  The map $f$ is  not an isomorphism in general. In fact two rays $\ell_1,\ell_2$ in $R(l)$ go to the same point $(\pi,x)\in \calP_l\times l$ if and only if $\pi = \la \ell_1,l) = \la\ell_2,l\ra$ and $\ell_1\cap \ell_2 = \{x\}$. This shows that the Schubert line 
$\sigma_{x,\pi}:=\sigma_x\cap \sigma_\pi$  is a secant line of $S$ passing through the point $x\in \bbG$. 

Let $r$ be the number of secant lines of $S$ in $\bbG$. It is called the \emph{rank} of $S$.  Note that, since $\bbG$  is a quadric hypersurface in $\bbP^5$, a secant line of $S$ through a general point in $\bbG$ is  contained in $\bbG$. Obviously $r = 0$ if $S$ is degenerate as a subvariety of $\bbP^5$. It is shown in \cite{Verra} that, if $S$ is nondegenerate, then $r = 0$ if and only if it is a Veronese surface in $\bbP^5$ realized as the congruence of lines of bidegree $(1,3)$, the congruence of secant lines of a twisted cubic.

It follows from the above that $R(l)$ is the normalization of a curve of  arithmetic genus $(m-1)(n-1)$ with $r$ double points. Therefore its genus $g$ is given by the  formula
\beq\label{rankformula}
g= (m-1)(n-1)-r.
\eeq

The composition of $f$ and the second projection is a degree $m$ cover of $\bbP^1$ of a curve of genus $g$ with $\deg(\Foc(S))$ ordinary branch points. Applying Hurwitz formula and \eqref{rankformula},  we obtain 
\begin{itemize}
\item[(2)] \it{The degree of the focal surface} is given by the formula
\beq\label{degreefocal}
\deg \Foc(S) = 2m+2g-2 = 2n(m-1)-2r,
\eeq
where $g$ is the sectional genus of $S$ and $r$ is the rank of $S$.
\end{itemize}

\textbf{2.2} Let $\sigma_{x,\pi} \subset \bbG$ be a secant line of $S$. The plane $\pi$  contains two rays of $S$ intersecting at the point $x$. It is called the \emph{null-plane} of the congruence and the point $x$ is its \emph{center}. Let $\textrm{Sec}(S) \subset \Gr_1(\bbP^5)$ be the variety of secant lines of $S$. The projection $\alpha:\textrm{Sec}(S)\to \bbP^3, \sigma_{x,\pi} \to x,$ is a map of degree $m(m-1)/2$, the projection $\beta:\textrm{Sec}(S)\to \check{\bbP}^3, \sigma_{x,\pi}\to \pi,$ is a map of degree $n(n-1)/2$:
$$\xymatrix{&\textrm{Sec}(S)\ar[dl]_\beta\ar[dr]^\alpha\\
\check{\bbP}^3&&\bbP^3.}$$
Note that, if  $m = 2$, the map $\alpha$ (resp. $\beta$) is an isomorphism over the set of nonsingular  points of $\Foc(S)$.

A general plane $\pi$ in the pencil $\calP_l$ contains $n(n-1)/2$ intersection points of $n$ rays in this plane. The locus of such points when $\pi$ moves in $\calP_l$ is a curve,  denoted in classical literature by  $|l|$. A plane $\pi$ containing $l$ contains $n(n-1)/2$ points equal to the pre-image of the corresponding point $\pi^*\in l^*$ under $\beta$ and also $r = \rank S$ points on $l$ that are the centers of null-planes containing $l$. In other words, the planes through $l$ define a pencil in $|l|$ with $n(n-1/2$ moving points and $r$ fixed point. Thus 
\beq\label{degreel}
\deg (|l|) = \frac{1}{2}n(n-1)+r.
\eeq

Let $x\in \Foc(S)$ be a singular point. We denote by $K(x)$ the cone of rays through $x$. Its degree is equal to $h(x)$.

Obviously,  $K(x)$ is an irreducible component of  the ruled surface $(l)$ when $x\in l$. The residual part is of degree $m+n-h(x)$. Let $C(x)$ be the curve in $S$ parameterizing the lines in the ruling of $K(x)$. Its degree  in the Pl\"ucker embedding is equal to $h(x)$. Since $R(l)$ is a hyperplane section of $S$, it intersects $C(x)$ with multiplicity $h(x)$. Thus $(l)$ and $K(x)$ intersect in $S$ with multiplicity $h(x)$. If $h(x) > 1$, then each pair of common rays lies in the same plane $\la l,x\ra$, hence $x\in |l|$.  This shows that $x$ is a $\half h(x)(h(x)-1)$-multiple point of $|l|$. It is also a point of multiplicity $h(x)$ on $(l)$.

Since each point of $|l|$ lies on a ray intersecting $l$, the  curve $|l|$ is contained in the ruled surface $(l)$. Under the normalization map $q'{}^{-1}(C(l))\to (l)$ the pre-image of a general point on $|l|$ consists of $\half n(n-1)$ points.

\begin{itemize}
\item[(3)] \it{The curve $|l|$  passes through  singular points $x$ such that  $h(x) > 1$ with multiplicity $\half h(x)(h(x)-1)$. It is contained in the locus  of singular points of multiplicity $\half n(n-1)$ on $(l)$.}
\end{itemize}

Next we take a plane $P^*$ dual to a point $P\in \bbP^3$ and define the surface $(P) = \alpha(\beta^{-1}(P^*))$ in $\bbP^3$. It follows from the definition, that the surface $(P)$ is the locus of points on some line through $P$ which are centers of null-planes containing the line. The duality argument of the proof of formula \eqref{degreel} shows that 
\beq\label{degP}
\deg (P) = \half m(m-1)+r.
\eeq

For any  singular point $x$ of $\Foc(S)$ with $h(x) = 1$, the plane-cone $K(x)$ does not contain a general point $P$. Thus the surface $(P)$ does not contain $x$\footnote{Jessop on p. 262 of \cite{Jessop} apparently mistakenly says that $(P)$ passes through all singular points.}. On the other hand, if $h(x) > 1$, the plane $\pi$ spanned by $P$ and a ray $\ell$ passing through $x$ contains $h(x)-1$ additional rays $\ell'$ through $x$, hence $x$ lies in $(P)$. The null planes spanned by $\ell$ and the remaining $n-h(x)$ rays in $\pi$ have centers belonging to $(P)$. Since a general $\ell$ intersects $(P)$ at $\deg(P)$ points, we obtain the following. 

\begin{itemize}
\item[(4)] \textit{Each singular point $x\in \Foc(S)$  is a singular point of $(P)$ of multiplicity 
$\half m(m-1)+r-n+h(x).$}
\end{itemize}

\textbf{2.3} Let $B$ be the branch curve of the degree 2 map $q':R_{S'}\to S'$. Let $\deg B = h\cdot [B]$ be the  degree of the image of $B$ on $S$ in the Pl\"ucker embedding. Using \eqref{omegaZ} and \eqref{ramdiv} and applying the adjunction formula, we get 
\beq\label{canrs}
K_{R_{S'}} = (K_{Z_{S'}}+R_{S'})\cdot R_{S'} = q'{}^*(2K_{S'}+2h).
\eeq
 The Hurwitz type formula now gives
$2K_{R_{S'}} = q'{}^*(2K_{S'}+B)$, hence 
\beq\label{branchcurve}
[B] =  2K_{S'}+4h.
\eeq

A ray $\ell_s$ may be tangent to $\Foc(S)$ with multiplicity 4 at a nonsingular point. It is clear that in this case $q'{}^{-1}(\phi^{-1}(s))$ intersects $R_{S'}$ at one point, hence $s' = \phi^{-1}(s)$ belongs to $B$. Conversely, if $s'$ belongs to $B$, then the ray $\ell_{\phi(s')}$ intersects $\Foc(S)$ at one point with multiplicity $4$.  Let $\Foc(S)_0$ be the set of such points in $\Foc(S)$. The pre-image of $\Foc(S)_0$ under the projection $p':R_{S'}\to \Foc(S)$ is equal to the ramification locus of $q':R_{S'}\to S'$. Since 
$q'{}_*p'{}^*\calO_{\bbP^3}(1) = \calO_{S'}(1)$, the projection formula shows that the degree of $\Foc(S)_0$ is equal to  $\deg B$. Formula \eqref{branchcurve} gives 
$$\deg B = 2K_{S'}\cdot h+4h^2 = 2(K_{S'}+h)\cdot h+ 2h^2 = 2(2g-2)+2(m+n).$$
Taking into account formula \eqref{rankformula}, we obtain

\begin{itemize}
\item[(5)]  \textit{The degree of the locus of points $\Foc(S)_0$  on $\Foc(S)$  where a ray touches  $S$ with multiplicity $4$ is equal to 
$\deg B$}.  
\end{itemize}

For any singular point $x$ of $\Foc(S)$, the fiber $q^{-1}(x)$ is a connected curve equal to the intersection  $C(x) = \sigma_x\cap S$.  Thus  $K(x)$ is  a cone over $C(x)$ of degree $h(x)$. Its pre-image under $p:Z_S\to \bbP^3$ is equal to the ruled surface $ q^{-1}(C(x))$ with base curve $C(x)\subset S$.  Its pre-image in $Z_{S'}$ is the ruled surface $q'{}^{-1}(\tilde{C}(l))$, where $\tilde{C}(l)\subset S'$ is the pre-image of $C(l)$ under the normalization map $\phi:S'\to S$. Assuming that  a singular point of $\Foc(S)$ is an ordinary node, the exceptional curve $E_x = p'{}^{-1}(x) \subset R_{S'}$  must be a smooth rational curve and, under the projection $q':Z_{S'}\to S'$, it is mapped isomorphically onto $\tilde{C}(x)$ and defines the normalization map $E_x \cong \tilde{C}(x) \to C(x)$. Since $C(x)$ is a rational plane curve of degree $h(x)$, it is singular if $h(x) > 2$, and hence $\phi:S'\to S$ is not an isomorphism over singular points of $C(x)$.
 
\begin{itemize} 
\item[(6)] \textit{Suppose $S$ is  smooth and  $x$ is an ordinary node of $\Foc(S)$. Then $h(x)\le 2$}.
\end{itemize}

In fact, under the generality assumption that all singular points of the curve $C(x) = \sigma_x\cap S$ are ordinary nodes, we obtain that the number  of singular points of $S$ is equal to 
\beq\label{numbersingular}
s= \half \sum_{x\in \Foc(S)}(h(x)-1)(h(x)-2)
\eeq

\section{Congruences of degree $2$ and class $n$}\label{S:congruences(2,n)}
\textbf{3.1} We specialize assuming that $m = 2, n\ge 2$.  We also assume that the congruence does not have fundamental curves. This excludes only one case when $(m,n) = (2,2)$ and $S$ is a quartic scroll  (see \cite{ArrondoGross}). We have in this case $r = 1, g = 0$ and $\deg\Foc(S) = 2$. So may assume that $\deg\Foc(S) \ge 4$.

\begin{itemize}
\item[(7)] \it{$\deg\Foc(S) = 4$ and $ g = 1, r = n-2$}.  
\end{itemize}

Let us see that $\deg\Foc(S) = 4$, the rest follows from  formula \eqref{degreefocal}. Suppose the degree is larger. Then a general  ray $\ell$ is tangent to the focal surface at two points, hence  $\ell$ intersects it at some other (nonsingular) point $x$. Since $m = 2$ and $x$ belongs to $\Foc(S)$, there will be another ray $\ell'$ tangent to $\Foc(S)$ at the point  $x$.  Thus there are two different rays through $x$ contradicting the definition of the focal surface.

Since $g = 1$, a general hyperplane section $H$ of $S$ is an elliptic curve. By the adjunction formula, the anti-canonical sheaf $\omega_{S'}^{-1}\cong \calO_{S'}(1)$ is ample. By definition, $S'$ is a del Pezzo surface. The anti-canonical linear system $|-K_{S'}|$ embeds $S'$ into $\bbP^d$, where 
$d = K_{S'}^2$ is equal to the degree of the image. It follows from (5) that $B\in |-2K_{S'}|$ and $\deg B = 2K_{S'}^2 = 2(2+n)$, hence we obtain the following.

\begin{itemize}
\item[(8)]\textit{The normalization $S'$ of $S$ is a del Pezzo surface of degree $d = 2+n$. The normalization map $\phi:S'\to S$ is an isomorphism if $n = 2,3$ and is the projection of the surface $S'$ anti-canonically embedded in $\bbP^{2+n}$ to a surface of degree $2+n$ in $\bbG \subset \bbP^5$. The branch curve $B$ of $q':R_{S'}\to S'$ belongs to $|-2K_{S'}|$.}
\end{itemize}
Collecting the previous formulas, we obtain the following.
\begin{itemize}
 \item[(9)] {The degree of the curve $|l|$ is equal to $\half (n^2+n-4)$. The degree of the  surface $(P)$ is equal to $n-1$}
 \item[(10)] \it{The surface $(P)$ is the locus of centers of null-planes that contain the point $P$. It passes through singular points $x$ of $\Foc(S)$  with multiplicity $h(x)-1$}.
\end{itemize}

Since the degree $d$ of a del Pezzo surface in $\bbP^d$ satisfies $d\le 9$, we see that $n\le 7$.  Moreover, there are two different types of congruences of class $6$. They correspond to two types of del Pezzo surfaces of degree $8$: re-embedded quadrics and blow-ups of one point in $\bbP^2$.

We know from \eqref{canrs} that 
$\omega_{R_{S'}} \cong \calO_{R_{S'}}$. Since $\Foc(S)$ is a quartic surface in $\bbP^3$, we have $\omega_{\Foc(S)} \cong \calO_{\Foc(S)}$, hence $p'{}^*\omega_{\Foc(S)} = \omega_{R_{S'}}$, hence all singular points of the quartic surface $\Foc(S)$ are rational double points. Since $Z_{S'}$ is smooth, all singular points of $\Foc(S)$ are singular points of the congruence $S$. Also,  no ray $\ell$ is contained in $\Foc(S)$ (otherwise, each $x\in \ell$ is contained in two rays, $\ell$ and another ray from $C(\ell)$ passing through $x$). This implies that the map $q:R_{S'}\to S'$ is a finite morphism of degree 2, and hence the branch curve $B$ is smooth.

\begin{itemize}
\item[(11)] \it{The surfaces  $R_{S'}$ is a K3 surface. All singular points of $\Foc(S)$ are double rational points  and the map $p':R_{S'}\to \Foc(S)$ is a minimal resolution of singularities. The double cover $q':R_{S'}\to S'$ is a finite morphism of degree $2$ with smooth branch curve $B\in |-2K_{S'}|$.}
\end{itemize}

Since a ray $\ell_s$ is tangent to $\Foc(S)$ at two points, $p^{-1}(\ell)$ splits into the union of two smooth rational curves in $Z_S$, one of them is the fiber $q^{-1}(s)$. The other component parametrizes rays that intersect $\ell_s$. It is projected to a hyperplane section of $S$ with a singular point at $x$. This curve is cut out by a hyperplane tangent to $\bbG$   at the point $s\in S$.

From now on we assume that $\Foc(S)$ has only ordinary nodes as singularities. 

Two ruled surfaces $(l)$ and $(l')$ intersect at $n+m = n+2$ common rays. Thus the two surfaces residually intersect along a curve of degree $(n+2)^2-(n+2) = n^2+3n+2$. If $x$ is a point on the intersection curve not lying on the common rays, then it is the intersection point of two different rays, hence no other ray passes through it. This shows that three ruled surfaces can intersect only at singular points of $\Foc(S)$ and at $3(n+2)^2$ points on the common rays (a ray common to two must meet the third). Since each singular point $x\in \Foc(S)$  has multiplicity $h(x)$ on $(l)$, the intersection curve has multiplicity $h(x)^2$ at $x$. Intersecting with the third ruled surface $(l'')$, we get 
$$(n+2)^3 = 3(n+2)^2+\sum i^3\alpha_i,$$
where $\alpha_i$ is the number of singular points $x$ with $h(x) = i$.
The Table \ref{table1} of possible solutions can be found in \cite[p. 280]{Jessop}.

\begin{table}[h!]
\centering
\begin{tabular}{|c|c|c|c|c|c|c|c|}\hline 
&$(2,2)$&$(2,3)$&$(2,4)$&$(2,5)$&$(2,6)_I$&$(2,6)_{II}$&$(2,7)$\\ \hline
$\alpha_1$&16&10&6&3&1&0&0\\
$\alpha_2$&&5&6&6&4&8&0\\
$\alpha_3$&&&2&3&6&0&10\\
$\alpha_4$&&&&1&0&4&0\\
$\alpha_5$&&&&&1&&\\
$\alpha_6$&&&&&&&1\\ \hline
$\sum \alpha_i$&16&15&14&13&12&12&11\\ \hline
\end{tabular}
\vskip 5pt
\caption{Singular points of congruences of bidegree $(2,n)$}\label{table1}
\end{table}

Applying formula \eqref{numbersingular}, we see that all congruences $S$ are singular if $n > 3$, and the number of singular points is equal to $\half (n-2)(n-3)$.

\textbf{3.3} Let $x$ be a singular point of $\Foc(S)$. Then  a general ray of the cone $K(x)$ is tangent to $\Foc(S)$ at some point outside $x$. The closure of the locus of the tangency points is a \emph{trope-curve} $T(x)$ on $\Foc(S)$  of degree $2h(x)$.  It is cut out in $\Foc(S)$  with multiplicity 2  by the cone $K(x)$. If $h(x) = 1$ (resp. $h(x) = 2$), it is a conic (resp. quartic curve) called  a \emph{trope-conic} (resp. \emph{trope-quartic}) of $\Foc(S)$. Since by (11) the projection $R_{S'}\to S'$ is a finite morphism, no ray lies in $\Foc(S)$. This implies  the  the trope-curves are irreducible.

The orbits of the birational involution of $\Foc(S)$ corresponding to the deck transformation $\sigma$ of the cover $q':R_{S'}\to S'$ are pairs of tangency points of a ray of the congruence. If $x\in \textrm{Sing}(\Foc(S))$, then the exceptional curve $E_x\subset R_{S'}$ is mapped under the involution to some other curve $\sigma(E_x)$. The curve $E_x+\sigma(E_x)$ is equal to the pre-image of the curve $\tilde{C}(x)= \phi^{-1}(C(x))\subset S'$ under the covering map $q'$. This shows that the curve $\tilde{C}(x)$ splits under the cover, and hence  is everywhere tangent to the branch curve $B$.
The degree of $\tilde{C}(x)$ in the anti-canonical embedding of $S'$ is equal to $h(x)$.

   Two singular points of $\Foc(S)$ are called \emph{conjugate} if the line joining these points is a ray from the congruence. Suppose $x,x'$ are two conjugate singular points.  Then the ray $\ell=\la x,x'\ra$ is contained in the intersection of the cones $K(x)$ and $K(x')$ and joins their vertices. A general plane $\pi$ containing $\ell$ intersects $K(x)$ along $h(x)$ rays and intersects $K(x')$ along $h(x')$ rays. Since $\pi$ contains $n$ rays, we see that the points $x$ and $x'$ are conjugate if $h(x)+h(x') > n$. In this case the ray $\la x,x'\ra$ must have multiplicity $h(x)+h(x')-n$.  Note that the inequality $h(x)+h(x')> n$ is only a sufficient condition for the conjugacy. 

\begin{itemize}
\item[(12)] \it{Suppose $h(x)+h(x') > n$, then the points $x,x'$ are conjugate and the ray joining the two points has multiplicity $h(x)+h(x')-n$.} \end{itemize} 

The \emph{conjugacy graph} of $S$ is a graph whose vertices are singular points of $\Foc(S)$ and two conjugate vertices $x,x'$ are joined by an edge. We also mark a vertex $x$ with the integer equal to $h(x)$ (if its larger than $1$) and mark the edge with $(h(x)+h(x')-n$ if $h(x)+h(x')-n> 1$.

\begin{itemize} 
\item[(13)] \it{Let $\tilde{C}(x)$ be the pre-image in $S'$  of the curve $C(x)$ on $S$, where $x$ is a singular point of $\Foc(S)$. Then $\tilde{C}(x)$ is everywhere tangent to the branch curve $B$ and has degree $h(x)$ in the anti-canonical embedding of $S'$. The conjugacy graph of $S$ is the dual intersection graph of the curves $\tilde{C}(x)$ on $S'$. }
\end{itemize}

\begin{example} Assume $n = 4$. The surface $\Foc(S)$ has two singular points $x_1,x_2$ with $h(x_i) = 3$. Thus $C(x_1)$ and $C(x_2)$ are  irreducible singular cubic curves and $S$ is singular with one singular point corresponding to the multi-ray that passes through the points $x_1$ and $x_2$. It is a common line of the cones $K(x_1)$ and $K(x_2)$. The normalization $S'$ of $S$ is a del Pezzo surface of degree $6$ anti-canonically embedded in $\bbP^6$. The normalization map $\phi:S'\to S$ is the projection from a general point on the secant hypersurface of $S'$. The projection from a general point of $\bbP^6$ is a smooth del Pezzo surface of degree $6$ that is contained in a smooth quadric. If we realize this quadric as the Klein quadric $\bbG$, we obtain a realization of a smooth del Pezzo surface of degree $6$ as a smooth congruence of lines if bidegree $(3,3)$.
\end{example}

 \section{Congruences of type $(2,3)$}\label{congruences(2,3)}
\textbf{4.1} Let us now summarize what we have found in the case $n = 3$. We assume that $S$ is general enough so that all singularities of the focal surface are ordinary nodes. 
\begin{itemize}
\item[(14)] \it{The surface $S$ is a smooth quintic del Pezzo surface anti-canonically embedded in $\bbG\subset \bbP^5$}.
\item[(15)] \it{The focal surface $\Foc(S)$ is a quartic surface with 15 nodes, ten of them with $h(x) = 1$ and remaining ones with $h(x) = 2$}. 
\item[(16)] \it{The surfaces $(l)$ are quintic elliptic ruled surfaces. A curve  $|l|$ is a singular curve of multiplicity $3$  on $(l)$, it coincides with the base curve of the pencil of quadrics generated by $(P)$ and $(P')$.} 
\item[(17)] \it{The surfaces $(P), P\in \bbP^3,$ form a web $W$ of quadrics  with base points at the set of nodes $x$ with $h(x) = 2.$} Its discriminant surface $\sfD_W$ parameterizing singular quadrics coincides with $\Foc(S)$.
\item[(18)] \it{The surface $Y = R_{S'}$ is a K3 surface. The morphism $p'=p:Y\to \Foc(S)$ is a minimal resolution of singularities and the morphism $q' = q:Y\to S'=S$ is a finite morphism of degree $2$ with a smooth branch curve $B\in |-2K_S|$. }
\end{itemize}

Let us denote the set of nodes $x$ of $\Foc(S)$ with $h(x) = 1$ (resp. $h(x) = 2$)  by $\calL$ (resp. $\calC$).  The curves $C(x),x\in \calL,$ are lines on $S$  tangent to $B$ at one point. The curves $C(x),x\in \calC,$ are conics in $S$ tangent to $B$ at two points.  The curves $T(x) = (K(x)\cap \Foc(S))_{\textrm{red}}$ are trope-conics if $x\in \calL$ and trope-quartics if $x\in \calC$. 

Since $\sigma(E_x)\cdot E_x = h(x)$, we see that $x$ is a simple point of $T(x)$ if $x\in \calL$ and a double point if $x\in \calC$. We know that two singular points $x,x'$ with $h(x) = h(x') = 2$ are conjugate, i.e. the cones $K(x)$ and $K(x')$ have a common ray. This implies that  $T(x)$ and $T(x')$ intersect and the corresponding conics $C(x)$ and $C(x')$ belong to different pencils of conics on $S$. 
 
 Let us recall  that a quintic del Pezzo surface $S$ is isomorphic to the blow-up of four points $p_1,p_2,p_3,p_4$ in the plane, in general linear position. Let $e_0$ be the class of the pre-image of a line and $e_1,\ldots,e_4$ be the classes of the exceptional curves  over the points $p_i$. The anti-canonical embedding is given by the linear system $|3e_0-e_1-e_2-e_3-e_4|$ represented in the plane by the linear system of cubic curves passing through the points $p_1,\ldots,p_4$. In the anti-canonical embedding,  the surface has  $10$ lines with the classes $e_1,\ldots,e_4$ and $e_0-e_i-e_j$. It has $5$ pencils of conics with classes $e_0-e_i$ and $2e_0-e_1-e_2-e_3-e_4$. 
 
  We denote by $\phi_{\calL\calL}:\calL\to 2^\calL$ the map that assigns to a point $x\in \calL$ the set of points from $\calL$ contained in $\sigma(E_x)\setminus \{x\}$. We have similar maps $\phi_{\calC\calC}:\calC\to 2^\calC, \phi_{\calL\calC}:\calL\to 2^\calC$ and $\phi_{\calC\calL}:\calC\to 2^\calL$. It follows from the above   that 
$$\#\phi_{\calL\calL}(x) = 3, \quad \#\phi_{\calC\calL}(x) = 4, \quad \#\phi_{\calL\calC}(x) = 2, \quad \#\phi_{\calC\calL}(x) = 4.$$

\begin{figure}[ht]
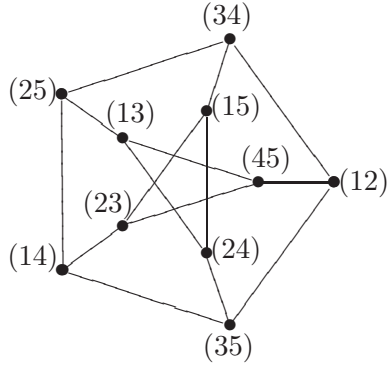
\label{peterseng}
\centering
\xy (-60,0)*{};(20,0)*{\bullet};(6.2,19)*{\bullet}**\dir{-};(-16.2, 11.7)*{\bullet}**\dir{-};(-16.2, 11.7)*{\bullet};(-16,-11.7)*{\bullet}**\dir{-};
(-16.2,-11.7)*{\bullet};(6.2,-19)*{\bullet}**\dir{-};(6.2,-19)*{};(20,0)*{\bullet}**\dir{-};
(10,0)*{\bullet};(-8.1,5.85)*{\bullet}**\dir{-};(-8.1,5.85)*{};(3.1, -9.5)*{\bullet}**\dir{-};(3.1, -9.5)*{};(3.1,9.5)*{\bullet}**\dir{-};
(3.1,9.5)*{};(-8.1,-5.85)*{\bullet}**\dir{-};
(-8.5,-5.85)*{};(10,0)*{}**\dir{-};   (10,0)*{};(20,0)*{}**\dir{-};  (3.1,9.5)*{};(6.2,19)*{}**\dir{-};
 (-16.2, 11.7)*{};(-8.1,5.85)*{}**\dir{-};(-16.2, -11.7)*{};(-8.2,-5.85)*{}**\dir{-};(6.2, -19)*{};(3.1,-9.85)*{}**\dir{-};
 (24,0)*{(12)};(5.8,22)*{(34)};(-20,12)*{(25)};(6,-22)*{(35)};(11,3)*{(45)};(7,10)*{(15)};(7,-9)*{(24)};(-7,9)*{(13)};
 (-10,-3)*{(23)};(-20,-10)*{(14)};
\endxy
\caption{The Petersen graph}
\label{petersen}
\end{figure}

In the sequel we will be indexing the set  $\calC$ by the set  $\{16,26,36,46,56\}$ and the set $\calL$ by subsets of $\{1,2,3,4,5\}$ of cardinality 2 as in the Petersen graph \ref{petersen}. We find that 
\begin{eqnarray*}
\phi_{\calL\calL}(x_{ab}) &=& \{x_{cd}:(cd)\cap (ab) = \emptyset\},\\
\phi_{\calC\calC}(x_{a6}) &=& \{x_{b6}, a\ne b\}, \\
\phi_{\calL\calC}(x_{ab})&=& \{x_{a6},x_{b6}\},\\
\phi_{\calC\calL}(x_{a6})&=&\{x_{ab}, b\ne a,6\}.
\end{eqnarray*}

 \begin{itemize}
\item[(19)] \it{The conjugacy graph of the congruence $S$ is equal to the dual intersection graph of $10$ lines $l_{ab}$ on $S$ and $5$ conics $C_i$, each taken from one of the five pencils of conics. It is equal to the union of the Petersen graph and the complete graph $K(5)$ on the set $\{1,2,3,4,5\}$. Each vertex $(ab)$ of the Petersen graph is connected to two vertices $a,b$  of $K(5)$ and each vertex $a$ of $K(5)$ is connected to 4 vertices $(ab),b\ne a,$ of the Petersen graph}.\end{itemize}

\textbf{4.2} The plane $K(x), x\in \calL,$ that cuts out a trope-conic passes through $x$ and the nodes from $\phi_\calL(x)$ and $\phi_{\calL\calC}(x)$. This shows that, for any $x\in \calL$,  

\beq\label{hyperplane}
2\sigma(E_{x})\in |\eta_H -E_{x}-\sum_{x'\in \phi_{l}(x)}E_{x'}-\sum_{y\in \phi_{lc}(x)}E_y|.
\eeq
where $\eta_H = c_1(p^*\calO_{\bbP^3}(1))$.

Similarly, we find, that for any $E_y,y\in \calC,$
\beq\label{tropeconic}
2\sigma(E_{y})\in |2\eta_H -2E_y-\sum_{y'\in \calC,y'\ne y}E_{y'}-\sum_{x\in \phi_{cl}(y)}E_{x}|.
\eeq

Applying formulas \eqref{hyperplane} and \eqref{tropeconic}, we find that
\beq\label{involution}
\sigma(\eta_H) \in |4\eta_H-\sum_{x\in \calL}E_x-2\sum_{y\in \calC}E_y|.
\eeq
 The linear system $|\sigma(\eta_H)|$ is the restriction to $\Foc(S)$ of the 4-dimensional linear system of quartic surfaces with nodes on $\calC$ and  
simple points at $\calL$. By choosing an appropriate basis of the linear system, it defines a birational self-map of $\Foc(S)$ that  blows down $\sigma(E_x), x\in \calL\cup \calC$ to the  nodes.

Let $\tilde{B}\cong B$ be the ramification curve of $q:Y\to S$. We know that its image $\Foc(S)_0$ in $\Foc(S)$ is equal to the locus of points where a ray is touching $\Foc(S)$ with multiplicity 4. It follows from  (5) that it is a curve in $\bbP^3$ of degree $4(mn-r)-2(m+n) = 10$.  Since $\tilde{B}$ intersects each $E_{x}, x\in \calC,$ at two points and intersects the curves $E_{y}, y\in \calL,$ at one point, we see that it passes through all singular points of $\Foc(S)$. If we assume that the congruence is \emph{moduli general} in the sense that the  rational Picard group $\textrm{Pic}(S)_\bbQ$ is generated by the curves  $E_{x}, x \in \calL\cup \calC,$ and $\eta_H$ (see the explanation of this assumption in section 8), we obtain
\beq\label{tildeB}
\tilde{B} \sim \half (5\eta_H-2\sum_{x\in \calC}E_x-\sum_{x\in \calL}E_x) \sim \half (\sum_{x\in \calL}(E_x+\sigma(E_x)) \sim \half (\sum_{x\in \calC}(E_x+\sigma(E_x)).
\eeq
Applying formula \eqref{hyperplane}, we see that $\tilde{B}$ is a curve of degree 10 cut out with multiplicity 2 by a quintic surface with nodes at $y,y\in \calC,$ that passes through the remaining nodes with multiplicity 1.

The linear system $|\tilde{B}|$ defines a $\sigma$-equivariant map from $Y$ to a surface of degree 10 in $\bbP^6$. We will study this birational model in Section \ref{S:degree10}.

\textbf{4.3} Let $(P)$ be the quadric surface corresponding to a nonsingular point $P$ on $\Foc(S)$ and $\ell$ be the unique ray passing through $P$. Since a general point $x$ on $\ell$ intersects another ray, $\ell$ is contained in a null-plane. Since $P\in \ell$, we see that  $\ell$ is contained in $(P)$.  Conversely, if $l$ is a line on $(P)$, a general point on $\ell$ is the intersection point of two rays, hence $l$ lies on an infinite set of secants of $S$. This shows that $l$ is a ray. So, we obtain that the quadric surface $(P)$ contains only one line passing through $P$, hence it must be a quadric cone. Since the degree of the discriminant surface of a web of quadrics is equal to 4, we obtain that $\Foc(S)$ coincides with this surface.

 Also observe that, if $P$ is a singular point $P$ with $h(P) = 1$, then any ray through $P$ is contained in $(P)$, hence the plane $\la T(x)\ra = K(x)$ is contained in $(P)$. Thus the quadric $(P)$ is reducible. The singular line of the quadric $(P)$ is the tangent line to the trope-conic at the point $P$.

\begin{itemize}
\item[(20)] \it{The web $W$ of quadrics $(P), P\in \bbP^3,$ has the base locus equal to the set $\calC$ of singular points of $\Foc(S)$ with $h(x) = 2$. The surface $\Foc(S)$ coincides with the \emph{discriminant surface} $\sfD_W$ of $W$, i.e.  the locus of singular quadrics in $W$.  The 10 quadrics $(x)$, where $x\in \calL$, are reducible quadrics, one of the irreducible components is the plane  $K(x)$}.
\end{itemize} 

We will return to the discussion of the geometry related to the web of quadrics $W$ in Section \ref{S:15nodal}.

\section{The Segre cubic  and the Castelnuovo-Richmond-Igusa quartic}\label{segrecubic}
\textbf{5.1} Let us recall some known facts about the Segre cubic primal and the Castelnuovo-Richmond-Igusa quartic hypersurface that can be found, for example, in \cite[9.4]{CAG}. Let $|L|$ be the linear system of quadrics through the singular points $p_1,\ldots,p_5\in \bbP^3$ in general linear position. The linear system $|L|$ defines a map
\beq\label{kapranov}
\Phi:\bbP^3\da \bbP^4
\eeq
with the image  a cubic threefold in $\bbP^4$ projectively isomorphic to the  \emph{Segre cubic primal} $\sfS_3$ given by equations
$$\sum_{i=1}^6z_i^3 = \sum_{i=1}^6z_i = 0.$$ 
The permutation group $\frakS_6$ acts naturally on $\sfS_3$ via permuting the coordinates. We denote this action by ${}^I\frakS_6$. There is a natural isomorphism between $\sfS_3$ and the GIT-quotient of $(\bbP^1)^6$ by $\PGL(2)$. It follows from an explicit map from $(\bbP^1)^6$ to $\sfS_3$ given by the Joubert functions \cite[Theorem 9.4.10]{CAG} that the action of $\frakS_6$ on $\sfS_3$ via permuting the factors of $(\bbP^1)^6$ differs from the previous action by an outer automorphism of $\frakS_6$. We denote this action by ${}^{II}\frakS_6$. The linear action of $\frakS_5$ in $\bbP^3$ leaving invariant the set of points $p_1,\ldots,p_5$ embeds $\frakS_5$ onto a transitive subgroup of ${}^I\frakS_6$ and onto non-transitive subgroup of ${}^{II}\frakS_6$ (recall that there are two conjugacy classes of subgroups of $\frakS_6$ isomorphic to $\frakS_5$.)  

Recall that a \emph{duad}  is a subset of the set $[1,6] = \{1,\ldots,6\}$ of  cardinality $2$. A \emph{syntheme} is a set of three   complementary  duads.  A \emph{total} is a set of 5 synthemes that contain all duads. A subgroup of $\frakS_6$ fixing a total is isomorphic to $\frakS_5$ and acts transitively on the set $[1,6]$.

The images $P_{ij}$ of 10 lines $\la p_i,p_j\ra$ are singular points of $\sfS_3$, all ordinary nodes.  The images $\Pi_i$ of the exceptional divisors of the blow-up of $\bbP^3$ at the points $p_i$ and the planes $\la p_i,p_j,p_k\ra$ are planes contained in $\sfS_3$. The 15 planes and 10 nodes form an abstract configuration $(15_4,10_6)$. 

It is also known that the projectively dual hypersurface of $\sfS_3$ is projectively isomorphic  to the \emph{Castelnuovo-Richmond-Igusa quartic hypersurface} $\CR_4$. The equations of $\CR_4$ is 
$$4\sum_{i=1}^6y_i^4-(\sum_{i=1}^6y_i^2) = \sum_{i=1}^6y_i = 0.$$
The projective duality map 
\beq\label{dualitymap}
\Psi:\sfS_3\da \CR_4
\eeq
is given by the choice of a basis in the linear system of polar quadrics of $\sfS_3$. The  map $\Psi$ is equivariant with respect to the action  ${}^{II}\frakS_6$ on $\sfS_3$ and the action of $\frakS_6$ on $\CR_4$  via permuting the coordinates $y_i$. 

The images of planes in $\sfS_3$ under $\Psi$ are 15 double lines of $\CR_4$. Each line  intersects three other lines, and three lines pass through each of 15 intersection points forming  a symmetric abstract configuration $(15_3)$, the \emph{Cremona-Richmond configuration}. The images of the tangent cones at singular points of $\sfS_3$ are quadric surfaces in $\CR_4$ cut out by hyperplanes everywhere tangent to $\CR_4$ (\emph{cardinal hyperplanes}).

The Segre cubic $\sfS_3$ is characterized by the property that it contains 10 singular points, the maximal number for a normal cubic hypersurface in $\bbP^4$. The singular points are ordinary nodes forming a ${}^I\frakS_6$-orbit of the point  $[1,1,1,-1,-1,-1]$. They are indexed by subsets $(abc)$ of cardinality 3 of the set $[1,6]$, up to a complementary sets. The set of planes is the $\frakS_6$-orbit of the plane $\pi_{12,34,56}:z_1+z_2=z_3+z_4= 0$. They are indexed by synthemes $(ab,cd,ef)$. Two planes intersect along a line if and only if the synthemes share a duad.   Otherwise they intersect at a singular point of $\sfS_3$.

Dually, we have 15 double lines in $\CR_4$ forming a $\frakS_6$-orbit of the line given by equations $y_1-y_2 = y_3-y_4=y_5-y_6 = 0$. They are also indexed by synthemes. The intersection points of the double lines form a $\frakS_6$-orbit of the point $[-2,2,1,1,1,1]$. They are indexed by duads. Two double lines intersect if they share a duad.

A cardinal hyperplane is dual to a singular point of $\sfS_3$. The set of cardinal planes form one $\frakS_6$-orbit of the hyperplane $y_1+y_2+y_3-y_4-y_5-y_6 = 0$. They are indexed by subsets of $[1,6]$ of cardinality 3 up to complementary set.

The  cardinal hyperplane $H(123)$  intersects $\CR_4$ along the quadric given by equation $y_1y_2+y_1y_3+y_2y_3-y_4y_5-y_4y_6-y_5y_6 = 0$. A cardinal hyperplane $H(abc)$ contains 6 double lines corresponding to the synthemes $(ai)(bj)(ck)$, where $\{i,j,k\}$ is the complementary set to the set $\{a,b,c\}$.

 We will denote by $\bbP_{\textrm{Seg}}^4$ the projective space $|L|^*$, where $\sfS_3$ lies, and denote  by $\bbP_{\CR}^4$  the dual space $|L|$, where $\CR_4$ lies,   
 
 A general  hyperplane $H$ in $\bbP_{\CR}^4$ cuts out $\CR_4$ along a 15-nodal quartic surface $X_{\CR}^H$. Its pre-image $X_{\Seg}^H$ under $\Psi$ is cut out in $\sfS_3$ by a quadric $Q^H$. It is the polar quadric of $\sfS_3$ with the pole  $h$ corresponding to $H$ via the projective duality.

The pre-image of $X_{\Seg}^H=\sfS_3\cap Q^H$ under $\Phi$ from \eqref{kapranov} is a quartic surface $K^H$ with nodes at the points $p_1,\ldots,p_5$. The linear system of such surfaces defines the composition map 
$$\Psi\circ \Phi:\bbP^3\da \bbP_{\CR}^4.$$ 
The surface $K^H$ contains two sets of 10 lines, ten of them are the lines $\la p_i,p_j\ra$ and the other ten are the residual lines in the sections of $K^H$ by the planes $\la p_i,p_j,p_k\ra$. In a minimal resolution of $K^H$ we get two sets of disjoint $(-2)$-curves, each one from one set intersects three members of the second set.  

\begin{itemize}
\item[(21)] \textit{A $15$-nodal quartic surface $X_{\CR}^H$ is birationally isomorphic to a quartic surface $K^H$ in $\bbP^3$ with $5$ nodes and two sets of disjoint lines forming a symmetric configuration 
$(10_3)$}.
\end{itemize}

\textbf{5.3}
It is known that the Fano variety $F(\sfS_3)$ of lines in $\sfS_3$ consists of $21$ irreducible components. Fifteen of them parameterize lines in 15 planes. The remaining six components  are isomorphic to the quintic del Pezzo surface $\dP_5$ \cite{DolgSegre}. We denote them by 
$F(\sfS_3)_1,\ldots,F(\sfS_3)_6$. The six components are transitively permuted by ${}^{II}\frakS_6$. The stabilizer subgroup of $F(\sfS_3)_i$ in ${}^{II}\frakS_6$ is isomorphic to $\frakS_5 = \Aut(\dP_5)$.  The five del Pezzo components $F(\sfS_3)_i, i = 1,\ldots,5,$ parameterize lines in $\sfS_3$ that meet the planes $\Pi_i$. Under the map  $\Phi:\bbP^3\da \sfS_3$ from \eqref{kapranov}, they are realized as the blow-ups of $4$ points on the plane $\Pi_i$  equal to the images of lines $\la p_i,p_j\ra, j\ne i$. The component $F(\sfS_3)_6$ parameterizes lines that meet five planes $\pi_1,\ldots,p_5$. By \emph{Theorem on the fifth associated plane}  (see \cite[Chapter X, Theorem XXIV]{Semple}, or \cite[Theorem 8.5.3]{CAG}) this is a quintic del Pezzo surface.

Thus the choice of a non-transitive subgroup $\frakS_5$ of ${}^{II}\frakS_6$ or a transitive subgroup of ${}^{I}\frakS_6$ defines the choice of $5$ planes among the $15$ planes contained in $\sfS_3$. As we noted earlier, this choice is equivalent to the choice of a total. 

\begin{itemize}
\item[(22)] \textit{The choice of a total in the set $[1,6]$ defines the choice of $5$ planes in $\sfS_3$, or, equivalently, $5$ double lines in $\CR_4$, or equivalently, $5$ nodes of a $15$-nodal quartic $X_{\CR}^H$. There are six such choices.}
\end{itemize}

The intersection of the quadric $Q^H$ with a general line on $\sfS_3$ corresponding to a point in $F(\sfS_3)_i$ consists of two points. This defines a birational involution on $X_{\Seg}^H$ (or $X_{\CR}^H$) and hence a biregular involution on the K3-surface $Y^H$ birationally isomorphic to $X_{\Seg}^H$ and $X_{\CR}^H$. Its orbits correspond to points in $F(\sfS_3)_i$.

\begin{itemize}
\item[(23)] \textit{The choice of a total defines a birational involution on $X_{\CR}^H$ and a biregular involution on $Y^H$ with quotient isomorphic to the quintic del Pezzo surface $\dP_5$}.
\end{itemize}

\textbf{5.4} Let us recall some known facts about the geometric interpretation of irreducible representations of $\frakS_5$ (see, for example, \cite{Cheltsov}, \cite{DFL}). Let $D$ be the union of $10$ lines on a quintic del Pezzo surface $\dP_5$ and $\Omega_{\dP_5}(\log D)$ be the rank 2 vector bundle of logarithmic $1$-forms. It is invariant with respect to the group of automorphisms of $\dP_5$ isomorphic to $\frakS_5$. The linear space $H^0(\dP_5, \Omega_{\dP_5}(\log D))$ of global sections is an irreducible $5$-dimensional linear representation $W_5$ of $\frakS_5$. It generates the vector bundle and  defines a $\frakS_5$-equivariant anti-canonical embedding of $\dP_5$ into the Grassmannian $G(3,W_5) \cong G(2,W_5^\vee)$ by taking the kernels  (or cokernels) of the evaluation maps. The composition of this  map with the Pl\"ucker embedding defines an embedding $\dP_5\hookrightarrow \bbP(\bigwedge^2W_5^\vee)$ which is contained in the projectivization of a $6$-dimensional irreducible summand of $\bigwedge^2W_5^\vee$. 

The linear system 
$|\calO_{\bbP(W_5^\vee)}(2)| = \bbP(S^2(\bigwedge^2W_5))$ of quadrics in $\bbP(W_5^\vee)$ contains 
the subspace $|I_{\dP_5}(2)|$ of quadrics with base locus $\dP_5$. It corresponds to an irreducible summand of dimension $5$ of the linear representation $S^2(\bigwedge^2W_5)$ isomorphic to $W_5$. We will use that the representation $W_5$ of $\frakS_5$ is self-dual so that we may identify $G(2,W_5^\vee)$ with $G(2,W_5)$. The embedding $\dP_5\hookrightarrow G(2,W_5)$ is now defined by assigning to a point $x\in \dP_5$ the pencil of quadrics from $|I_{\dP_5}(2)|$ with singular point at $x$. Under the projection map of $\dP_5$ to $\bbP^4$ with center at $x$, this pencil is the pre-image of the pencil of quadrics containing the quartic del Pezzo surface equal to the image of the projection map (see \cite{Kapustka}).  

The universal family $Z_{\dP_5}$ of lines in $\bbP(W_5)$ parameterized by $\dP_5 \subset G(2,W_5)$ can be identified now with the incidence variety 
$$Z_{\dP_5} = \{(Q,x)\in |I_{\dP_5}(2)|\times \dP_5:x\in \textrm{Sing}(Q)\}$$ that comes with two projection maps 
$$\xymatrix{&Z_{\dP_5}\ar[dl]_{p_1}\ar[dr]^{p_2}& \\
\bbP(W_5) = |I_{\dP_5}(2)|&&\dP_5,}
$$
The image of the first projection $p_1$ consists of quadrics that have a singular point at some point in $\dP_5$. In fact, any singular quadric from $|I_{\dP_5}(2)|$ has this property. Otherwise, if we project from its singular point not in $\dP_5$, the image of a quintic surface $\dP_5$ will be contained in a quadric surface, obviously impossible (see \cite[Lemma 5.2]{Kapustka}). Thus, the image of $p_1$ is the  locus  $\Delta$ of singular quadrics in $|I_{\dP_5}(2)|$. The scheme-theoretical locus $\tilde{\Delta}$ of a linear system quadrics in $\bbP^5$ is given by vanishing of the discriminant polynomial of degree $6$. It is known that the  projective tangent space to a quadric  $Q$  of corank 1 considered as a point on  the  discriminant hypersurface of a linear system of quadrics is equal to the linear subspace of quadrics  that pass through the singular point of $Q$ \cite[(1.4.5)]{CAG}. Since all quadrics in $|I_{\dP_5}(2)|$ pass through $\dP_5$, we see that  a general quadric in $\Delta$ is a singular point of the discriminant scheme. This shows that 
$\tilde{\Delta} = 2\Delta$ and 
$$\deg \Delta = 3.$$ 
The restriction of the  logarithmic bundle $\Omega_{\dP_5}(\log D)$ to a line $l$ on $\dP_5$ is a rank two vector bundle of determinant $-K_{\dP_5}\cdot l = 1$. It follows that  the pre-image  $p^{-1}(l)$ is the minimal ruled surface $\bfF_1$ and the projection $p_1$ blows down its  exceptional section to a singular point of the image $\Delta$ of $p_1$. Thus we see that $\Delta$  is an irreducible cubic hypersurface in  $|I_{\dP_5}(2)|$  with 10 isolated singular points.  It must be isomorphic to the Segre cubic $\sfS_3$.
 
 \begin{itemize} 
 \item[(24)] \textit{There is a $\frakS_5$-equivariant isomorphism between the locus $\Delta$ of singular quadrics in $|I_{\dP_5}(2)|$ and the Segre cubic $\sfS_3$. The universal family of lines in $\sfS_3$ parameterized by a del Pezzo component of $F(\sfS_3)_i$ is $\frakS_5$-equivariantly isomorphic to the $\bbP^1$-bundle $Z_{\dP_5}\to \dP_5$.}
 \item [(25)]\textit{The projection $p_1:Z_{\dP_5}\to \sfS_3$ is a small projective resolution of the Segre cubic.}
 \end{itemize}

Our proof is a slight modification of the proofs from \cite[Lemma 2.31]{Cheltsov} and \cite[Theorem 5.6]{Kapustka}.

\begin{remark} In \cite{Finkelnberg} Finkelnberg proves that the Segre cubic admits six non-isomorphic small projective resolutions. Ours  represents one isomorphism class. It is one of the two isomorphism classes of small resolutions  that have the property that its automorphism group is isomorphic to $\frakS_5$.
\end{remark}

 We have already used the fact that   a quintic del Pezzo surface has 5 pencils $\calQ_i$ of conics. They define a conic fibration $\pi_i:S\to \bbP^1$. In the anti-canonical embedding the planes spanned by the members of the pencil $\calQ_i$ sweep a scroll $Z_i$. The scroll $Z_i$ is the image of the projective rank 2-bundle $\bbP(\calE_i)$, where $\calE_i = \pi_{i\ast}(\calO_S(-K_S))$. Applying the relative Riemann-Roch, we easily find that $c_1(\calE_i) \cong \calO_{\bbP^1}(1)^{\oplus 3}$, i.e. the projective bundle is isomorphic to $\bbP^2\times \bbP^1\to \bbP^1$. The map from this bundle to $\bbP^5$ is given by $|\calO_{\bbP(\calE_i)}(1)|$ 
and  coincides with the Segre map $ \bbP^2\times \bbP^1\to \bbP^5$.

Let $|\calO_{\bbP^5}(2)|\to |\calO_{Z_i}(2)|$ be the restriction map. Since $h^0(\calO_{Z_i}(2)) = h^0(S^2(\calE_i)) = h^0(\bbP^1,\calO_{\bbP^1}(2)^{\oplus 6}) = 18$, we obtain that each scroll $Z_i$ is contained in a net of quadrics. After we  identify the Segre variety $s(\bbP^2\times \bbP^1)$ with the projective space of $(2\times 3)$-matrices of rank 1, we easily see that the space of quadrics containing it consists of quadrics of rank 4. 

\begin{itemize}
\item[(26)]\cite{Kapustka} \textit{The space $|I_S(2)|$ of quadrics containing an anti-canonical quintic del Pezzo surface $S$ contains 5 planes $\Pi_i$ of quadrics of corank 2}.
\end{itemize}

 Using the projective duality, we can also embed  $\dP_5 \subset \Gr(2,W_5)$ in the dual Grassmannian $\Gr(3,W_5^\vee)$ of planes in the dual space $\bbP_{\CR}^4$. The universal family of planes parameterized by $\dP_5$ now becomes a $\bbP^2$-bundle over $\dP_5$. Its projection to $\bbP_{\CR}^4 = \bbP(W_5^\vee)$ is a double cover branched along the Castelnuovo-Richmond-Igusa quartic $\CR_4$ (see \cite[Proposition 2.4.5]{Cheltsov}).

\section{15 nodal quartic surfaces}\label{S:15nodal}
 
\textbf{6.1} It is a natural guess from reading the previous section that the choice of a total defines an isomorphism between the 15-nodal quartic surface $X_{\CR}^H$ and the focal surface $\Foc(S)$ of a congruence of lines $S$ of bidegree $(2,3)$. Note that this is not quite obvious since the construction gives a birational model of $X_{\CR}^H$ inside the universal family of lines in $\bbP^4$ (not in $\bbP^3$)  parameterized by $S$ considered as a surface of lines in $\Gr_1(4)$ (not in $\Gr_1(3)$).

Fix a component $F(\sfS_3)_i$ of $F(\sfS_3)$ and let $l$ be a line in it. We consider the plane $\la l,h\ra$ spanned by this line and the point $h\in \bbP_{\Seg}^4$ (dual to the hyperplane $H$).  The dual of this plane in $\bbP_{\CR}^4$ is a line $l^*$ contained in $H \subset \bbP_{\CR}^4$. The intersection $\sfS_3\cap \la l,h\ra$ consists of $l$ and a conic $C$ that intersect $l$ at two points. The intersection of the polar quadric $Q^H$  with  the plane $\la l,h\ra$ is equal to the intersection of the polar conic of $\sfS_3\cap \la l,h\ra$ with respect to the point $h$. It passes through its singular points, i.e. the intersection points of $l$ and $C$. Since the image of $Q^H$ under the duality map $\Psi$ is $H$, we obtain that the line $l^*$ intersects the quartic $X_{\CR}^H$ at two points, the images of $l\cap C$. In other words $l^*$ is a bitangent line of $X_{\CR}^H$. Thus we realize $F(\sfS_3)_i$ as an irreducible congruence of bitangent lines of the quartic surface $X_{\CR}^H$.

Take a general point $x\in H$. The set of lines $l^* \subset H$ containing $x$ corresponds to the set of planes $\la\ell,h\ra$ in the dual space that are contained in the hyperplane $x^\perp$ dual to the point $x$. This shows that the order $m$ of the congruence $F(\sfS_3)_i$ is equal to the number of lines in $\sfS_3$ that lie in  the  hyperplane $H$. If we consider the map \eqref{kapranov} and take $i = 6$ then this number is equal to the number of twisted cubic curves  parameterized by $F(\sfS_3)_6$ that lie in a quadric through the five points. This number is equal to $2$ and consists of the unique curves of bidegree $(1,2)$ and $(2,1)$   that pass through general $5$ points on the quadric. Similar count gives $m = 2$ for other components $F(\sfS_3)_i$ (in fact they all permuted by $\frakS_6$).

Formulas \eqref{rankformula} and \eqref{degreefocal} give now that the class $n$ of $F(\sfS_3)_i$ is equal to $3$. Alternatively, we can compute $n$  as follows. If we take a general plane $\pi$ contained in $H$, then its dual $\pi^\perp$ will be a line passing through $h$. It intersects $\sfS_3$ at three points. Through each of these points passes a unique line on $\sfS_3$ from the family $F(\sfS_3)_i$. These define three planes $\la l,h\ra$ and give three lines $l^*$ contained in $\pi$.

\begin{itemize}
\item[(29)] \textit{The choice of a total realizes $X_{\CR}^H$ as the focal surface of a congruence of bidegree $(2,3)$ in $\Gr_1(H)$}. 
\end{itemize}

 The surface $X_{\CR}^H$ also contains 10 tropes, the intersections with 10 cardinal hyperplanes. This gives ten more components of bi-degree $(0,1)$. It is known that the bidegree of the bitangent surface of a normal quartic surface is equal to $(12,28)$ \cite[p. 281]{Salmon}. Since $(12,28) = 10(0,1)+6(2,3)$ we have found all irreducible components.

\textbf{6.2} Let $X_{15}$ be any 15-nodal quartic surface in $\bbP^3$. Let us see how it arises as the  focal surface of a congruence of bidegree $(2,3)$ and hence it is isomorphic to a surface $X_{\CR}^H$ (another proof of this fact can be found in \cite{Avilov}). Projecting $X_{15}$ from a node, we obtain a double cover of $\bbP^2$ branched along a curve of degree 6 with 14 nodes. The only way to  realize this curve is to take the union of four lines $l_1,l_2,l_3,l_4$ in general linear position and a conic $C$ intersecting each line transversally outside the intersection points. The fifteenth node of $X_{15}$ arises, in the usual way, from a conic $K$ which everywhere tangent to the sextic curve. The linear system on the double cover $\pi:X_{15}\da \bbP^2$ that maps it to a quartic surface is equal to $|\pi^*(l)+K'|$ where $l$ is a line on $\bbP^2$ and $\pi^{-1}(K) = K'+K''$.

The reduced pre-images of the line components of the  branch curves are conics passing through the fifteenth node $P$. These are our trope-conics $\sigma(E_x), x\in \calL,$ passing through the node. In this way we see the configuration of 10 trope-conics and 15 nodes forming an abstract configuration of  type $(15_4,10_6)$.  The pre-image of the conic component $C$ of the branch locus is a curve of degree 4 passing doubly through $P$ and simply through 7 other nodes. Thus $X_{15}$ contains 15 nodes, 10 trope-conics passing through 6 nodes, and 5 quartic curves with a double point at one of the nodes. 

Consider the linear system   of polar cubics of the 15-nodal quartic $X_{15}$. On a minimal resolution $Y$ of $X_{15}$ it is defined by the linear system $|D| = |3\eta-\sum E_i|$, where $E_i$ are the exceptional curves of the resolution and $\eta$ is the pre-image of the class of a hyperplane section of $X_{15}$.  We check that its dimension is equal to 4 and its degree is equal to 6. 

We have $(D-E_i)^2 = 0,$ so that $|D-E_i|$  is an elliptic pencil and $E_i$ is mapped under the map $\phi_D:Y\to \bbP^4$ given by the linear system $|D|$ to a plane conic. The members of  $|D-E_i|$ are curves of arithmetic genus one and the map  $\phi_{|D|}$ maps these curves one-to-one to curves of arithmetic genus one spanning a $\bbP^3$ that contains the plane spanned by the conic image of $E_i$. This shows that the map $\phi_{|D|}$ is of degree $1$ and its image  $Z$ is a surface of degree $6$, the intersection  of a quadric $Q$ and a cubic  (not uniquely defined). The quadric $Q$ is tangent to a cubic at the points $q_1,\ldots,q_{10}$, the images of the $10$ trope-conics on $X_{15}$. Let us consider the 4-dimensional linear system of quadrics through $q_1,\ldots,q_{10}$. It defines a rational map $f:\bbP^4\da \bbP^4$. Its restriction to $Z$ is given by  the linear system $|A| = |2\eta_Z-R_1-\cdots-R_{10}|$ on the minimal resolution $\tilde{Z}$ of $Z$, where $\eta_Z$ is the pre-image of the class of a hyperplane section of $Z$ and $R_i$ are the exceptional curves over the points $q_i$. Since $A^2 = 4$, it maps $Z$ (and $\tilde{Z}$) to a quartic surface  in $\bbP^4$ contained in the hyperplane $H$ whose pre-image under the map is the quadric $Q$.  Thus the image $Z'$ of $Z$ is the intersection of a quartic hypersurface $V_4$ with the hyperplane $H'$. The preimage of $V_4$ under $f$ is a cubic hypersurface that contains $15$ planes spanned by the $15$ conics on $Z$. The intersection points of the $15$ planes give $10$ singular points of this hypersurface. By the projective uniqueness of the Segre cubic, the cubic hypersurface must be projectively equivalent to $\sfS_3$. Since $\CR_4$ is the image of $\sfS_3$ under the map $f$ given by quadrics through its nodes, we see that $f(\sfS_3)$ is projectively equivalent to $\CR_4$ and $Z'$ is projectively equivalent to  $X_{\CR}^H = \CR_4\cap H$ for some hyperplane $H$ projectively equivalent to $H'$. 

Let $\eta_{H}$ be the class of a hyperplane section of $X_{\CR}^H$. We have $\eta_Z = 3\eta-\sum E_i$ and 
$$\eta_H = 2\eta_Z-\sum R_i = 6\eta-2\sum E_i-\sum R_i$$
in $\Pic(Y)$. By expressing the classes $R_i$ corresponding to $10$ trope-conics on $X_{15}$ in terms of $\eta$ and the classes of $E_i$, and substituting into $\eta_H$ we obtain $\eta_H = \eta$. Thus we have constructed a projective isomorphism from $X_{15}$ to $X_{\CR}^H$.

\begin{itemize}
\item[(30)] \it{A 15-nodal quartic surface $X_{15}$ is  projectively isomorphic to a hyperplane section $X_{\CR}^H$ of the Castelnuovo-Richmond-Igusa quartic}. 
\item[(31)] \it{The congruence ${\rm Bit}(X_{15})$ of bitangent lines of $X_{15}$ consists of 16 irreducible components, six of them are congruences of bidegree $(2,3)$ and 10 are planes of bidegree $(0,1)$.}
\end{itemize}

\begin{remark} The last assertion should be compared with the similar assertion that the congruence of bitangent lines of a quartic Kummer surface with 16 nodes consists of 16 plane components and 6 components of bidegree $(2,2)$. The minimal resolution $R_S$ of the Kummer surface $\Foc(S)$ is a double cover of a quartic del Pezzo surface with branch curve $B\in |-2K_S|$. In the isomorphic model of $S$ as the blow-up of $5$ points $p_1,\ldots,p_5$ in the plane, the curve $B$ is the proper transform of a plane curve  of degree $6$ with cusps at the five points  and tangent to all lines $\la p_i,p_j\ra$ and the conic through the five points. The analog of the linear system  $|\tilde{B}|$ from \eqref{tildeB}  defines a smooth octic model of the focal surface $\Foc(S)$. It is equal to the intersection of three diagonal quadrics:
$$
\sum_{i=0}^5z_i^2 = \sum_{i=0}^5a_iz_i^2 = \sum_{i=0}^5a_i^2z_i^2 = 0,
$$
 the image of the  curve $\Foc(S)_0$ isomorphic to $B$ is a canonical model of $S$ obtained as the hyperplane section $z_i = 0$.  The six choices of the hyperplane correspond to the six realizations of the Kummer surface as the focal surface of a congruence of bidegeree $(2,2)$. The curve $B$ is a  \emph{Humbert curve} of genus 5 (see  \cite{Edge2}). It has a special property that it contains a maximal number (= 10) of vanishing even theta characteristics for a curve of genus 5 \cite{Varley}.  The number of moduli of Humbert curves is equal to 2.

\end{remark}

{\bf 6.3} Let $\calM_{4,k}$ be the moduli space of projective equivalence classes of $k$-nodal quartic surfaces in $\bbP^3$. It is a locally closed subset of the GIT-quotient of $|\calO_{\bbP^3}(4)|/\!/\textrm{SL}(4)$. It is an interesting problem to find out whether the irreducible components of $\calM_{4,k}$ are rational varieties.    It follows from the above  that any $15$-nodal quartic surface is birationally isomorphic to $X_{\CR}^H = \CR_4\cap H$. It is known that the group of projective automorphisms of $\sfS_3$ and $\CR_4$ is isomorphic to $\frakS_6$ (see for example \cite[\S 3]{Finkelnberg}). Thus  the stabilizer subgroup of $X_{\CR}^H$ in  $\Aut(\bbP_{\CR}^4) \cong \PGL(5)$ is isomorphic to $\frakS_6$. This gives  
\beq\label{moduli}
\bbP_{\textrm{Seg}}^4/\frakS_6 = \check{\bbP}_{\CR}^4/\frakS_6 \da \calM_{4,15}, \ H\mapsto H\cap \CR_4,
\eeq
admits the  inverse. It is known that 
$$\bbP_{\textrm{Seg}}^4/\frakS_6 \cong \bbP(2,3,4,5,6)$$
\cite[9.4.5]{CAG}. Therefore 
we obtain the following.

\begin{itemize}
\item[(32)] \it{The moduli space $\calM_{4,15}$ of 15-nodal quartic surfaces is a rational 4-fold}. 
\end{itemize}

Note that the moduli space of Kummer surfaces $\calM_{4,16}$ is realized as the hypersurface in $\calM_{4,15}$ birationally isomorphic to the quotient 
$$\CR_4/\frakS_6 \cong \bbP(2,3,5,6).$$
Of course, the rationality of $\calM_{4,16}$ is a well known fact.

\textbf{6.4}  Now we know that $X_{15}$ is projectively isomorphic to a quartic surface in $\bbP_{\CR}^4$ equal to a hyperplane section $X_{\CR}^H$ of $\CR_4$. So, we can identify the ambient space $\bbP^3$ with a hyperplane $H$ in $\bbP_{\CR}^4$. We can also identify $X_{15}$ with the focal surface $\Foc(S)$ of a congruence $S$ of lines in $H$ of bidegree $(2,3)$ that defines a partition $\calL+\calC$ of its $15$ nodes and defines the involution $\sigma$ such that $\sigma(E_x), x\in \calL,$ are its trope-conics and $\sigma(E_x),x\in \calC,$ are its trope-quartics. There are six ways to do it and we fix one of them. We will use the expressions of different linear systems on $X_{15}$ in terms of the divisor classes of $E_x,\sigma(E_x),$ on the minimal resolution $Y$ of $X_{15}$ given in subsection 4.2.

Recall from (17) that there is a web $W$ of quadrics $(P)$ in $H$ through the nodes $x\in \calC$ parameterized by points $P$ in $H$. So we can identify $W$ with $H$ and consider the map $f:X_{15}\da H^*$ given by the linear system $W$. It is a hyperplane in the linear system $|2\eta_H-\sum_{x\in \calC}E_x|$, where $\eta_H$ is a hyperplane section of $X_{15}$ (lifted to $Y$, as usual). We know from (17) that  $X_{15}$ coincides with the discriminant surface $\sfD_W$ of $W$. Recall that the \emph{Steinerian surface} (or the \emph{Jacobian surface}) of $W$ is the locus of singular points of quadrics from $W$ (see \cite[1.1.7]{CAG}). The rational map 
$$\st:\sfD_W \da \St(W)$$
that assigns to a singular quadric from $W$ its singular point is called the \emph{Steinerian map}. It is given by the linear system $|\eta_{\st}|$ such that 
$$2\eta_{\st} = 3\eta_H-\sum_{x\in \calL}E_x.$$
Since $\eta_{\st}^2 = 4$,  the Steinerian surface is a quartic surface in $H$. We see that it has $5$ nodes, the images of $E_x, x\in \calC,$ and $20$ lines, the images of $E_x, \sigma(E_x), x\in \calL$. Thus it coincides with the $5$-nodal quartic model $K^H$ of $X_{15}$ and we can identify the source space $\bbP^3$ in the map $\Phi$ from \eqref{kapranov} with the space $|\eta_{\st}|^*$. We have just described a \emph{Cremona transformation} of $H$ that transforms a $15$-nodal quartic to a $5$-nodal quartic. 

The map $\Phi$ is given now by the linear system  
\beq\label{steinerianmap}
|2\eta_{\st}-\sum_{x\in \calC}E_x|, 
\eeq
so  $W$ is identified with a hyperplane in this linear system. It maps  $X_{15}$ to the projection of the sextic surface $X_{\Seg}^H = \sfS_3\cap Q^H$ from the point $h\in \bbP_{\Seg}^4$ corresponding to $H$ via the projective duality. It is a sextic surface in the dual space $H^*$.

Next we shall identify this sextic surface with the dual surface $X_{15}^*$ of $X_{15}$. 

The pre-image of a plane in $H^*$ through a  point $x^*$ under the map $f$ is a quadric in $H$ that contains the fiber $f^{-1}(x^*)$. Thus this fiber is the base scheme of the net of quadrics in $H$ defined by the point $x^*$ considered as a plane in $H$. A general net of quadrics in $\bbP^3$ has $8$ distinct base points. In our case all quadrics pass through $5$ points from $\calC$. So, our map $f$ is of degree $3$. A net of quadrics has less than $8$ base points if one of the quadrics in the net has a singular point  at  a base point. This shows that the ramification divisor $\text{Ram}(f)$ coincides with the set of singular points of quadrics  from $W$. This is the Steinerian surface $\St(W)$.  On the other hand, a quadric is singular if it does not intersect transversally the branch divisor $f(\textrm{Ram}(f))$. Considered as a point in $H$, the quadric belongs to the discriminant surface $\sfD_W$ of $W$. We know from fact (17) that $\sfD_W$  coincides with our quartic surface $X_{15}$. Thus the branch divisor is the dual surface $X_{15}$* of $X_{15}$. The image of the Steinerian surface is the branch divisor. This gives an identification of the sextic surface in $H^*$ with the dual surface of $X_{15}$.

\begin{itemize}
\item[(33)] \it{If we identify a  $15$-nodal quartic surface $X_{15}$ with $X_{\CR}^H$, then its projectively dual surface $X_{15}^*$ is a surface of degree $6$ in $H^*$ obtained by projection of its sextic model $X_{\Seg}^H$ from the point $h$ dual to $H$}. 
\item[(34)] \it{The $5$-nodal quartic model $K^H$ is identified with the Steinerian quartic surface in $H$ of the web of quadrics $W$ whose discriminant surface is equal to $X_{15}$.} 
\end{itemize}
Note that degree of the dual surface  agrees with the Pl\"ucker formula \cite[Theorem 1.2.7]{CAG}, according to which the degree of the dual of a  $k$-nodal quartic surface is equal to $36-2k$. 

We refer to \cite[page 160]{Baker} for another more explicit map from a $X_{15}$ to its $5$-nodal quartic model.

\begin{remark} Note that $f^{-1}(f(\textrm{Ram}(f)) = 2\textrm{Ram}(f)+F$, and comparing the degrees, we obtain that the residual surface $F$ is of degree $10$. The map $f$ defines a birational isomorphism $F\cong X_{15}^*$, so it realizes a birational model of $X_{15}$ of degree $10$ in $\bbP^3$. In the next section we will study a nonsingular degree 10 model of $X_{15}$ in $\bbP^6$. The surface $F$ is its projection to $\bbP^3$. The surfaces $F$ and $\St(W) =\textrm{Ram}(f))$ intersect along a curve of degree $40$. We believe (but no proof yet) that this curve is birationally isomorphic to the curve $\Foc(S)_0$ when we identify $X_{15}$ with the focal surface of a congruence $S$ of bidegree $(2,3)$.
\end{remark}

\textbf{6.5} 
Recall  that a web $W$ of quadrics in $\bbP^3$ defines a \emph{Reye congruence} $\Rey(W)$ of lines in $\bbP^3$ (reducible if $W$ has base points). Its rays are  lines (called \emph{Reye lines}) which are contained in a subpencil of the web (see \cite{Cossec}, \cite[1.1.7]{CAG} and \cite[Chapter 8]{Enriques2}). If $W$ has no base points, then $\Rey(W)$ is irreducible and its  bidegree is equal to $(7,3)$. The lines through the base points are Reye lines, so, in our case, the Reye congruence has $5$ components of bidegree $(1,0)$. The remaining irreducible component $\Rey(W)_0$ is of bidegree $(2,3)$. 

Recall that a quadric $(P)\in W$ was defined in (17) as the locus of centers of null-planes containing $P$. Take a ray $l_s$ from $S$ that does not pass through the base points of $W$. Any general point $x\in \ell_s$ lies on another ray $\ell_{s'}$ such that $x = \ell_s\cap \ell_{s'}$ becomes the center of the null-plane $\la\ell_s,\ell_{s'}\ra$. This shows that the quadric $(x)$ contains the line $\ell_s$, and hence $\ell_s$ is a Reye line. This implies that  $\Rey(W)_0$  coincides with  $S$.  

A pencil $\calP_\ell$ of quadrics from $W$ containing a Reye line $\ell$ not through a base point has the base locus equal to the union of the  line and a twisted cubic curve intersecting it at two points and passing through the base points. We again see the congruence as a del Pezzo surface parameterizing twisted cubics through the set $\calC$.

The pencil of quadrics containing a general Reye line $\ell$ has exactly two singular quadrics $(P)$ and $(P')$ with singular points at the singular points of the base curve of the pencil. This shows that $\ell$ is tangent to the discriminant surface $X_{15}$ of $W$ at two points $P$ and $P'$. But it also intersects the Steinerian surface at two points. The first pair of points is an orbit of the involution $\sigma$ on $X_{15}$ defined by the congruence. The second pair of points is an orbit of an involution on $\St(W)$ which we denote by $\tau_{\rey}$ and call the \emph{Reye involution}. Note that they are different involutions on a nonsingular model $Y$ of both surfaces. As we know the former involution has the fixed locus an irreducible curve $\Foc(S)_0$ of degree $10$. The latter involution has fixed locus $5$ nodes of $\St(S)$ which become the disjoint union of $5$ smooth rational curves on $Y$. The base locus of the pencil in $W$ through a Reye line intersecting $\Foc(S)$ at a point in $\Foc(S)_0$ consists of this line and a twisted cubic tangent to the line. The base locus of the quadrics containing a Reye line passing through a base point $p_i$ consists of this line and a twisted cubic passing through all base points and tangent to the line at the point $p_i$. 

Finally observe that the proper transforms  $\sigma(E_x), x\in \calL,$ of ten trope-conics on  $Y$  are mapped under the Steinerian map to  lines joining two base points of $W$. They are Reye lines invariant under the Reye involution.

Let $Z =  Y/(\tau_{\rey})$ be the quotient surface by the involution $\tau_{\rey}$. Since its set of fixed points is the union of $5$ smooth rational curves $E_x,x\in \calC$, the usual formulas tell us that $Z$ is a smooth rational surface with 
Euler-Poncar\'e characteristic equal to $18$.  The branch curve of the quotient map $q:Y\to Z$ is the disjoint union of smooth rational curves $\bar{E}_x$ with self-intersection $-4$. The Hurwitz type formula $K_Y = q^*(K_Z)+\sum_{x\in \calC}E_x$ implies that 
$\sum_{x\in \calC}\bar{E}_x\in |-2K_Y|$. Since $h^0(-2K_Z) = 1$, we see that $h^0(-K_Z) = 0$. Thus $Z$ is a \emph{Coble surface}, a smooth projective rational surface with empty anti-canonical linear system and non-empty anti-bicanonical linear system. As we observed in the previous paragraph, the exceptional curves $\sigma(E_x), x\in \calL,$ are invariant under $\tau_{\rey}$. They intersect the ramification locus at two points from different irreducible components. Their images on $Z$ are $10$ disjoint $(-1)$-curves. After we blow down these curves we obtain a quintic del Pezzo surface $\bar{Z}$. The image of the branch locus is the union of $5$ conics on it (in an anti-canonical model). Since we identified Reye lines with rays of the congruence $S$, the surface $\bar{Z}$ can be identified with $S$ and the conics with the conics $C(y), y\in \calC,$ corresponding to trope-conics on $\Foc(S) = X_{15}$. Fix a realization of $X_{15}$ as the focal surface $\Foc(S)$ of a quintic del Pezzo congruence of lines in $\bbP^3$. 

\begin{itemize}
 \item[(35)] \it{The blow up of $10$ intersection points of $5$ conics $C(y), h(y) = 2,$  on $S$ is a Coble surface $Z$ with $K_Z^2 = -5$. The double cover of $Z$ branched at the proper transforms of the conics is isomorphic to a K3 model $Y$ of $X_{15}$. The covering involution coincides with the Reye involution $\tau_{\rey}$ associated to the web $W$ of quadrics $(P), P\in \bbP^3,$ with discriminant surface $X_{15}$.}
\end{itemize} 
 
 Recall that any quintic del Pezzo surface $S$ contains $5$ pencils of conics. If we choose one conic $C_i$ from each pencil and blow up the ten intersection points we obtain a Coble surface $Z$ with $|-2K_Z|$ represented by the union of proper transforms $\bar{C}_i$ of the conics. The double cover branched along this divisor is a K3 surface with $5$ isolated $(-2)$-curves, the reduced pre-images of the curves $\bar{C}_i$. There are also $10$ disjoint $(-2)$-curves, the pre-images of the exceptional curves over the intersection points. Note that the construction depends on $5$ parameters defined by a choice of a conic from each pencil. Note that degenerating $5$ conics to the union of $5$ reducible conics, we obtain a `most algebraic' K3 surface with Picard number $20$ studied by E. Vinberg \cite{Vinberg}. The isomorphism classes of our  K3 surfaces depend on $4$ parameters, so they define  a hypersurface in $(\bbP^1)^5/\frakS_5 \cong \bbP^5$. What is its degree?
 
 \begin{remark} It is known that in the case when a web $W$ of quadrics has no base points and its Steinerian surface is nonsingular, the Reye congruence  $\Rey(W)$ is an Enriques surface realized as a congruence of lines of bidegree $(7,3)$ and the discriminant surface $\sfD_W$ is a Cayley quartic symmetroid with $10$ nodes. The map  
$$\nu:\Rey(W)\to \Bit(\sfD_W), $$ 
from $\Rey(W)$ to the surface $\Bit(\sfD_W)$ of bitangents of $\sfD_W$ that assigns to a Reye line $\ell$ the pencil of quadrics $\calP_\ell$ in $W$ containing $\ell$ is the normalization map \cite[8.4]{Enriques2}.  In our case, $\Bit(\sfD_W)$ is a reducible surface withe each irreducible surface isomorphic to a quintic del Pezzo surface. If we fix one corresponding to $S$ from the above, then the map  $Z\to S$ described in (35) is the map $\nu$ over an irreducible component $S$ of 
$\Bit(\Foc(S)) = \Bit(\sfD_W)$.
\end{remark}

\begin{remark} The fact that a $15$-nodal quartic surface  is a quartic \emph{symmetroid} (i.e. equal to the determinant of a symmetric matrix with linear forms as its entries) should be compared with the well known fact that the Kummer surface $X_{15}$ is also a quartic symmetroid with the Steinerian surface  equal to the $6$-nodal \emph{Weddle quartic surface}. 
However in the Kummer case the web of quadrics is not formed by surfaces $(P)$ because it follows from formula \eqref{degP} that in this case they are planes but not quadrics. 

Let us explain the relationship between the Weddle and the Kummer surface. Recall that we have a map $\Phi:\bbP^3\da \sfS_3$ given by the linear system of quadrics through $5$ points $p_1,\ldots,p_5$. If we add one more point $p_6$, then the web $W$ of quadrics through the six points defines a map
$\phi_W:\bbP^3 \da W^*$ 
and its image is equal to the projection of the image $\sfS_3$ of the map $\Phi:\bbP^3\da \sfS_3$ from the point $\Psi(p_6)\in \sfS_3$. As we explained in subsection 6.4, the Steinerian surface $\St(W)$ is the ramification divisor of the map $\phi_W$. It is a quartic surface with $6$ nodes at $p_1,\ldots,p_6$, classically called a \emph{Weddle surface}. The branch divisor of the map is the dual surface of the discriminant surface $\sfD_W$ of $W$.  

On the other hand, we know that, the image of the ramification divisor under the map $\Psi$ is the intersection $Q\cap \sfS_3$, where $Q$ is the polar quadric with the pole at $x  = \Phi(p_6)$ that corresponds to the hyperplane $H(y)$ that is tangent to $\CR_4$ at the point $y = \Psi(x)$ and cuts out our Kummer surface $X_{16}$. This identifies $W^*$ with $H(y)$ and the branch divisor with $X_{16}$, hence the discriminant surface $\sfD_W$ with $X_{16}^*$. Since $X_{16}$ and $X_{16}^*$ are projectively isomorphic, we get a symmetric determinantal realization of any Kummer surface. There are many classical sources that explain this construction (see, for example, \cite[page   124]{Coble}).

We have described the birational map from the Weddle surface $\St(W)$ to the Kummer surface $X_{16}$ as the map from the ramification divisor of $\phi_W$ to the branch divisor of the degree two map given by the web $W$ quadrics through the nodes of $\St(W)$. The inverse map from $X_{16}$ to $\St(W)$ is given by cubic surfaces that pass through $10$ nodes of $X_{16}$ complementary to a set of nodes on one of the trope-conic that contains the node $y$ of $X_{16}$ (a \emph{odd trope-conic}).   It maps the six nodes to the points $p_1,\ldots,p_6$ which are the six nodes of the Weddle surface $\St(W)$. It maps the trope-conic to the twisted cubic  through the points $p_1,\ldots,p_6$ and it maps all other $15$ trope-conics to lines $\la p_i,p_j\ra$ (see \cite[Chapter XI, \S 97]{Hudson}).  The Reye involution is now the composition 
$$\tau_{\Rey} = \sigma\circ p_y\circ \sigma,$$ 
where $\sigma$ is the switch involution corresponding to one of the six realizations of $X_{16}$ as the focal surface of a congruence of lines of bidegree $(2,2)$ and $p_y$ is the projection involution of $X_{16}$ from the point $y$.
The locus of fixed points of the  Reye involution is the union of $6$ tropes passing through the node $y$. The quotient by the Reye involution is a Coble surface obtained by blowing up $15$ intersection points of six lines in general position in the plane. The surface is birationally isomorphic to the double cover of the plane branched along the six lines. This is of course a familiar birational model of a Kummer surface (see \cite[Theorem 10.3.16]{CAG}).

Note that the Kummer surface admits another symmetric determinant realizations  as the Hessian surface of a nonsingular cubic surface (see \cite[Theorem 4.1]{DK}) which has no  analog for a $15$-nodal quartic surface. 
\end{remark}

\textbf{6.5} Let us introduce some elliptic pencils on $Y$ invariant with respect to the Reye involution (recall that we have six of them corresponding to a choice of an isomorphism $X_{15}\cong \Foc(S)$).

Let
\beq\label{reyeroot}
\frakr  = 2\eta_H-\sum_{x\in \calL}E_x = 2\eta_{\st}-\eta_H \in \Pic(Y),
\eeq
where we use the definition of $\eta_{\st}$ from \eqref{steinerianmap}. It defines a reflection involution of $\Pic(Y)$
\beq\label{reyeroot2}
s_{\frakr}: v\mapsto v+\frac{x\cdot \frakr}{2}\frakr.
\eeq
One checks that $\tau_{\Rey}^* = s_{\frakr}$ \cite[10.4]{Enriques2}. In particular, we have 
\begin{eqnarray*}
\tau_{\Rey}^*(E_x) &\sim&2\eta_H-\sum_{x'\in \calL,x'\ne x}E_x \quad \textrm{if $x\in \calL$},\\
\tau_{\Rey}^*(E_y) &=& E_y \quad \textrm{if $y\in \calC$},\\
\tau_{\Rey}^*(\sigma(E_x)) &= &\sigma(E_x)\quad \textrm{if $x\in \calL$},\\
\tau_{\Rey}(\sigma(E_y))&\sim &\sigma(E_y)+4\eta_H-2\sum_{x\in \calL}E_x \quad  \textrm{if $y\in \calC$},\\
\tau_{\Rey}^*(\eta_H)&=& 8\eta_{\st}-3\eta_H,\\
\tau_{\Rey}^*(\eta_{\st}) &=& 3\eta_{\st}-\eta_H.
\end{eqnarray*}

For any line $l$ on $\St(W)$, the pencil of plane sections of $\St(W)$ containing $l$ defines the residual pencil of cubic curves. It gives rise to  an elliptic pencil 
$ |\eta_{\st}-L|$, where $L$ is the proper inverse transform of $l$ on $Y$. Thus we have 10 elliptic pencils $|F_x| = |\eta_{\st}-E_x|, x\in \calL$ and 10 pencils
$|F_x'| = |\eta_{\st}-\sigma(E_x)|$.  Using the formulas from the above, we verify that 
$$\tau_{\Rey}^*(F_x) = \tau_{\Rey}^*(\eta_{\st}-E_x) = \eta_{\st}-E_x.$$
This shows that the pencils $|F_x|$ are invariant with respect to the Reye involution. We have
\beq
\sum_{x\in \calL}F_x = 10\eta_{\st}-\sum_{x\in \calL}E_x = 10\eta_{\st}+2\eta_{\st}-3\eta_H = 3(4\eta_{\st}-\eta_H),
\eeq
The linear system
\beq\label{fanosystem}
\big|\frac{1}{3}\sum_{x\in \calL}F_x\big |= |4\eta_{\st}-\eta_H|
\eeq
maps $Y$ to a surface of degree 20 in $\bbP^{11}$. The ramification divisor of  the projection $Y\to Z$ is equal to $R = \sum_{y\in \calC}E_y$. Since $F_x\cdot E_y = 0$, each pencil $|F_x|$ descends to an elliptic pencil $\bar{F}_x$ on the Coble surface $Z$. The linear system \eqref{fanosystem}  descends to a linear system 
$$|\Delta| = \big|\frac{1}{3}\sum_{x\in \calL}\bar{F}_x\big|$$
that maps $Z$ to a rational surface of degree 10 in $\bbP^5$, a \emph{Fano model} of the Coble surface $Z$ (see \cite[Chapter 10]{Enriques2}). It has $5$ ordinary nodes, the images of the curves $E_y,y\in \calC$.

 Note that the involution $\sigma$ defined by the congruence of liens $S$ transforms $|F_x|$ to the pencil $|F_x'|$. It is invariant with respect to the conjugate involution $\sigma\circ\tau_{\rey}\circ \sigma$.

Let us now consider other elliptic pencils invariant with respect to the Reye involutions. They are the pre-images $|F_y|, y\in \calC,$ of the pencils of conics on $S$ under the composition map $f:Y\to Z = Y/(\tau_{\rey}) \to S$. Recall that each pencil contains the curve $C(y), y\in \calC,$ defined by the trope-quartic on $X_{15} = \Foc(S)$.

Take one conic $C(y), y\in \calC,$ from each of the five pencils. Since the map ramifies over $C(y)$, its pre-image  on $Y$ is the curve $E_y, y\in \calC$ that enters in the fiber with  multiplicity $2$.  The other components are the pre-images of four curves $\sigma(E_x), x\in \calL,$ that intersect $\sigma(E_y)$ with multiplicity one. This shows that the fiber is of type $\tilde{D}_4$ (or type $I_0^*$ in Kodaira's notation). In our duad notations of the curves $E_x,x\in \calL\cup \calC$, we re-denote $|F_y|$ with $|F_a|$, where $y =(a6)$. Using formula from subsection 4.2, we find that 
\beq\label{newpencil}
F_a \sim \half(4\eta_H-\sum_{b\ne 6} E_{y}-\sum_{b\ne a} E_{b6}-2\sum_{c,d\ne a,6} E_{cd})
\eeq
and $F_a\cdot F_b = 2$. Adding up, we get a linear system
\beq\label{deg10map1}
|\half(F_1+\cdots+F_5)| = |5\eta_H-\sum_{y\in \calC}E_{y}-2\sum_{x\in \calL}E_{x})| = 
 |\sum_{x\in \calL}\sigma(E_x)+\sum_{y\in \calC}E_y)|\eeq
which defines a $\tau_{\Rey}$-equivariant birational morphism to  a degree 10 surface in $\bbP^6$. It blows down  the curves $\sigma(E_x),x\in \calL,$ to $10$ singular points, maps the curves $E_x, x\in \calL,$ to rational quartic curves, maps the curves $E_y, y\in \calC,$ to conics and maps the curves  $\sigma(E_y), y\in \calC,$ to rational curves of degree $6$.

\begin{itemize}
\item[(36)]\it{The surface $X_{15}$ admits 15 elliptic pencils $|F_x|,x\in \calL\cup \calC,$ invariant with respect to the Reye involution $\tau_{\rey}$. The linear system 
$|\frac{1}{3}\sum_{x\in \calL} F_x|$ descends to a linear system on the Coble surface $Z = Y/(\tau_{\rey})$ that defines a Fano model of $Z$ as a rational $5$-nodal surface of degree $10$ in $\bbP^5$. The linear system $|\frac{1}{2}\sum_{y\in \calC} F_y|$ defines a $10$-nodal degree $10$ birational model of $Y$ in $\bbP^6$.}
\end{itemize}

\section{Degree 10 model of a 15-nodal quartic surface}\label{S:degree10}
\textbf{7.1} We continue to identify a 15-nodal quartic surface $X_{15}$ with the focal surface $\Foc(S)$ of a quintic del Pezzo congruence of lines $S$.  We have two  linear systems $|\tilde{B}|$ given by equation \eqref{tildeB} and given by equation \eqref{deg10map1}. The former linear system is $\sigma$-invariant and ample. It maps the  minimal resolution $Y$ of $X_{15}$ onto a nonsingular K3 surface of degree 10 in $\bbP^6$. The latter is $\tau_{\rey}$-invariant and maps $Y$ onto a $10$-nodal degree $10$ surface in $\bbP^6$. In this section we will only consider the nonsingular degree 10 model $Y_{10}$ of $Y$ obtained in this way.

 In the anti-canonical embedding of $S$ the curve $B\in |-2K_S|$ is cut out in $S$ by a quadric. Since $S$ is a linear intersection of the Grassmannian $\Gr_1(4)$, the curve $B$ is a canonical curve in $\bbP^5$ equal to the intersection of quadrics. By Petri Theorem, it  is neither hyperelliptic,  nor trigonal, no contains a $g_2^5$. It also embeds in $S$ by global sections of the logarithmic bundle on $S$. It follows from \cite{mukai} that it has only finitely many $g_4^1$, i.e. it is not a bielliptic curve.  By another result of Mukai \cite{Mukai},  the image of $Y$ under the map defined by $|B|$ is the intersection of the Fano threefold $\sfV_5$ of genus 6 (degree 5) and index 2 in $\bbP^6$ with a quadric $Q$. The images of all curves $E_x,\sigma(E_x), x\in \calL$ are lines and the images of $E_x,\sigma(E_x), x\in \calL,$ are conics.

We know that, for any $y\in \calC$, we have $E_y\cdot \sigma(E_y) = 2$. This means that $G_y = E_y+\sigma(E_y)$ varies in an elliptic pencil $|G_x|$ on $Y$. The curves $E_x,E_{x'},\sigma(E_x),\sigma(E_{x'}), x,x'\in \calL$ such that 
$\sigma(E_x)\cdot E_y = \sigma(E_{x'})\cdot E_y = 0$ form a reducible fiber of $|G_y|$ of type $\tilde{A}_3$. There are 3 such fibers. The curves $E_{x}, x\in \calL,$ such that $E_x\cdot E_y = 1$ define 4 sections of the elliptic fibration.

We have $G_y\cdot G_{y'} = 2$ and it follows from \eqref{tildeB} that 
$$2\tilde{B} \sim \sum_{y\in \calC}G_y.$$
The elliptic pencils descend to pencils of conics on the quotient $S$ by the involution $\sigma$. It follows that the curve $\tilde{B}$ passes through the intersection points $E_i\cap E_i'$ and two opposite singular points in each reducible fiber of the pencil.
Since $\tilde{B}\cdot G_y = 4$, the image of each member  of $|G_y|$ in the degree 10 model $\sfV_5\cap Q \subset \bbP^6$ of $X_{15}$  is a quartic curve of arithmetic genus one, and as such it spans a $\bbP^3$ in $\bbP^6$. Since $G_y\cdot G_{y'} = 2$, we see that any pair of the $\bbP^3$'s intersects along a subspace of dimension $\ge 1$.

\textbf{7.2} Consider the scroll $\calS_i$ swept by the 3-planes spanned by the images of the members of the elliptic pencil $|G_y|$. It is the image of the projective bundle $\bbP(\calE_i)$, where $\calE_i = (f_i)_*(\calO_Y(\tilde{B}))$ is a rank 4 vector bundle on $\bbP^1$. We have already observed in section 5.4  that a  quintic del Pezzo surface $S$ equipped with a conic bundle structure $\pi_i:S\to \bbP^1$ embeds in the trivial projective bundle $\bbP^1\times \bbP^2 = \bbP(\calN_i)\to \bbP^1$, where $\calN_i = \pi_{i\ast}(\omega_S^{-1}) \cong \calO_{\bbP^1}(1)^{\oplus 3}$. The quotient projections $q:Y\to S$ show that 
$$\calE_i \cong \pi_{i\ast}(q_*(\calO_{Y}(\tilde{B}))) = \pi_{i\ast}(\omega_S^{-1}\oplus \calO_S) = 
\calO_{\bbP^1}(1)^{\oplus 3}\oplus \calO_{\bbP^1}.$$
The scroll $\calS_i$  is equal to the image of $\bbP(\calE_i^\vee)$ under the map given by the linear system $|\calO_{\bbP(\calE)}(1)|$ (see also \cite{Knutsen} for the confirmation of this fact.\footnote{Note that there is a misprint in the inequality for $c= \textrm{Cliff}(C) < [\frac{g-1}{2}|$ for the Clifford index for a curve $C$ of genus $g$ in \cite[Theorem 3.2]{Knutsen}. According to \cite{Green}, it should be $\textrm{Cliff}(C) \le [\frac{g-1}{2}]$. In our case $C = \tilde{B}$ is a curve of genus 6, and we get $c = 2$.}) 

The image of the section $\bbP^1\to \bbP(\calE_i)$ defined by the projection $\calE_i\to \calO_{\bbP^1}$ is blown down to a singular point $s_i$ of $\calS_i$ which is contained in all rulings of the scroll. The projections $\bbP^6\da \bbP^5$ with center at $s_i$ maps the scroll to the scroll $\bbP(\calN)$. Its restriction to $Y_{10}$ is the projection $q:Y_{10}\to S$. 

The canonical sheaf of $\bbP=\bbP(\calE^\vee)$ is equal to $\calO_{\bbP}(-4)\otimes \pi_i^*\omega_{\bbP^1}$, where $\pi_i$ is the projection to the base. Using the adjunction formula, we see that $Y_{10}\subset \calS_i$ is equal to the complete intersection of the images of two divisors $D_1,D_2\in |\calO_{\bbP}(2)\otimes \pi_i^*\calO_{\bbP^1}(-1)|$ under the projection to $\bbP^6$ (see \cite[9.2]{Knutsen}). The linear system $|\calO_{\bbP}(2)\otimes \pi_i^*\calO_{\bbP^1}(-1)|$ is cut out by quadrics in $\bbP^6$ that contain a ruling 3-plane of $\calS_i$. Since a smooth quadric in $\bbP^6$ does not contain a 3-plane, all quadrics that cut out $Y_{10}$ in $\calS_i$ must be singular. They all pass through the singular points $s_i$ of $\calS_i$.

\begin{itemize}
\item[(37)]\it{The nonsingular $\sigma$-invariant degree 10 model $Y_{10}$ of a 15-nodal quartic $Y$ is contained in five scrolls $\calS_i$ swept by 3-planes spanning a fiber of  the elliptic fibrations $f_i:Y\to \bbP^1$ given by the pencil $|G_y|, y\in \calC$. Each scroll is the image of the projective bundle $\bbP(\calE) = \bbP(\calO_{\bbP^1}(1)^{\oplus 3}\oplus \calO_{\bbP^1})$ in $\bbP^6$ under the linear system $|\calO_{\bbP(\calE)}(1)|$. The section defined by the surjection $\calE\to \calO_{\bbP^1}$ is blown down to a singular points $s_i$ of $\calS_i$. The surface $Y_{10}$ is cut out in $\calS_i$ by two quadrics singular at the point $s_i$.}
\end{itemize}

\textbf{7.3}  It is known that $\sfV_5$ is equal to a transversal intersection of the Grassmannian $G = \Gr_1(\bbP(V_5))$ of lines in $\bbP^4 = |V_5|$ with a linear subspace $|W_7|$ of dimension 6 in the Pl\"ucker space $|\bigwedge^2V_5|$. The linear system $|I_G(2)|$ of quadrics containing $G$ is of dimension 4 and consists of quadrics of rank 6. It can be identified with the projective space $|V_5|$ via assigning to a vector $v\in V_5$ the quadratic form $q(x) = v\wedge x\wedge x$. After we choose a volume form on $V_5$, this corresponds to an isomorphism $I_{G}(2) \cong \bigwedge^4V_5^\vee \cong V_5$. 

Let $Y$ be a transversal intersection of $\sfV_5$ with a quadric $Q = V(q)$. It is a smooth K3 surface of degree 10 and genus 6 in $\bbP^6= |W_7|$.  We have 
\beq\label{decompose}
I_Y(2) = I_{\sfV_5}(2)\oplus \la q\ra
\eeq  Since a general quadric in  $I_{\sfV_5}(2)$ is of corank 1, we obtain that the discriminant hypersurface of degree 7 of singular quadrics in $|I_Y(2)|$ is equal to the union of a hyperplane $\calH = |I_{\sfV_5}(2)|$ and a hypersurface $\calD_Y$ of degree 6.

Now let $Y = Y_{10}$. Consider the restriction of the linear system of quadrics in $\bbP^6$ to a scroll $\calS_i$. We have 
$$S^2(\calE) = S^2(\calO_{\bbP^1}(1)^{\oplus 3})\oplus S^2(\calO_{\bbP^1})\oplus \calO_{\bbP^1}(1)^{\oplus 3}.$$ 
This implies that 
$\dim H^0(\bbP(\calE),\calO_{\bbP(\calE)}(2)) = 6h^0(\calO_{\bbP^1}(2))+3h^0(\calO_{\bbP^1}(1)+1 = 25.$ This shows that $\dim |I_{\calS_i}(2)| = 28-25 = 3$, hence each scroll $\calS_i$ is contained in a net of quadrics. Since each quadric containing $\calS_i$ contains a pencil of 3-planes and also contains $Y_{10}$, we obtain that the discriminant hypersurface of $\bbP^5 = |I_{Y_{10}}(2)|$.  In fact, projecting a quadric $\calQ$ containing $\calS_i$ from the singular point $s_i$ of $\calS_i$, we find that the fibers contain $\bbP^3$ as a codimension one subvariety, hence the corank of $\calQ$ must be greater than or equal to  $2$.

\begin{itemize}  
\item[(38)] \it{The discriminant hypersurface of the linear system $|I_{Y_{10}}(2)|$ of quadrics containing $Y_{10}$ contains 5 planes in its singular locus parameterizing quadrics of corank $\ge 2$}. 
\end{itemize}

 This should be compared with an analogous fact for a quintic del Pezzo surface $S$ discussed in subsection (7.2). 
 
\textbf{7.4} Let $Y$ be again any transversal intersection $\sfV_5\cap Q$. Using \eqref{decompose}, we obtain a decomposition
$$ \bigwedge^3I_Y(2) = \bigwedge^3I_{\sfV_5}(2)\oplus \bigwedge^2I_{\sfV_5}(2)\otimes \la q\ra = 
\bigwedge^3V_5\oplus \bigwedge^2V_5\otimes \la q\ra.
$$ 
If we view $G = \Gr_1(\bbP^4)$ as $\Gr(2,V_5)$ for some linear space of dimension $5$, then we have a natural (means $\SL(V_5)$-equivariant) isomorphism 
$$V_5 \to I_G(2), \quad v\mapsto v\wedge x\wedge x.$$
 We view the quadric $Q = V(q)$ as a restriction of a quadric on the Pl\"ucker space $\bigwedge^2V_5$, so, via polarization with respect to a vector $v\in  \bigwedge^2V_5$, it defines a linear form $L_v\in \bigwedge^2V_5^\vee \cong \bigwedge^3V_5$. We associate to $Q$ the linear subspace
$$A_Q:= \{(L_v,v\otimes q), v\in \bigwedge^2V_5\} \subset  \bigwedge^3I_Y(2).$$
One checks that $A_Q$ is a 10-dimensional Lagrangian subspace of $\bigwedge^3I_Y(2)$ with respect the wedge-product pairing, and as such defines an \textit{\textrm{EPW}-sextic hypersurface} in $|I_Y(2)|$ isomorphic to the sextic hypersurface $\calD_Y$ \cite[Proposition 2.1 and its proof]{Iliev}.

It is known that the irreducible family of EPW-sextics coming from a general degree 10 K3 surfaces via the construction from the above form a codimension 1 subvariety in the moduli space of EPW-sextics. 
The dual EPW-sextics of those coming from K3's have a singular point of multiplicity 3 \cite[Proposition 3.4]{O'Grady}.

 Let $Y^{[2]}$ be the \emph{Hilbert scheme} of $0$-cycles $\xi$ on $Y$ of length 2. 
For any such cycle $\xi$, its span $\la \xi\ra$ in $\bbP^6$ is a line. The linear subspace  of quadrics in $|I_Y(2)|$ containing this line is a hyperplane. This defines a rational map
\beq\label{alpha}
\alpha:Y^{[2]}\da |I_Y(2)|^*.
\eeq
It is clear that it is not defined on the set of cycles $\xi$ that are contained in a line in $Y$. This set is the union of planes $P_l = l^{[2]}$, where $l$ is a line in $Y$.  Obviously  any line on $Y$ is contained in $\sfV_5$. It is known that the Hilbert scheme $\Fano_1(\sfV_5)$ of lines  in $\sfV_5$ is isomorphic to $\bbP^2$ \cite[Corollary (6.6)]{Iskovskikh}. Each  line $l\in \Fano_1(\sfV_5)$ either intersects $Y$ at a $0$-cycle $\xi\in Y^{[2]}$ or is contained in $Y$. Thus $\Fano_1(\sfV_5)^0 = \Fano_1(\sfV_5)\setminus \{\textrm{lines in}\ Y\}$ embeds in $Y^{[2]}$. Let 
$\widetilde{Y^{[2]}}$ be the blow-up of the planes $P_l$ and $\tilde{P}$ be the closure in $\widetilde{Y^{[2]}}$ of $\Fano_1(\sfV_5)^0$.

It is known that the Hilbert scheme  of conics $\textrm{Fano}_2(\sfV_5)$ is isomorphic to $\bbP^4$ \cite{Iliev0}.  By \cite[Lemma 4.20]{O'Grady0}, for any $\xi\in Y^{[2]}$ not contained in a line in $\sfV_5$, there is a unique conic $C_\xi$ on $\sfV_5$ that contains $\xi$. The intersection of $C_\xi\cap Y$ is cut out by the quadric $Q$, thus one can define the residual set $0$-cycle $\xi'$ of length 2 (O'Grady shows that the residual $0$-cycle is well defined for a non-reduced $Z$). This defines a birational involution $\iota$ of $Y^{[2]}$.

Since $C_\xi$ is contained in $\sfV_5$, the set of quadrics in $I_Y(2)$ containing the plane $\la  C_\xi\ra$ is a hyperplane. Thus $\alpha(\xi) =\alpha(\iota(\xi))$ is the same point in $W$. This shows that the degree of $\alpha$ is at least 2.  

The following results can be found in \cite{O'Grady}. 

\begin{itemize} 
\item[(39)] \it{For a general $Y = \sfV_5\cap Q$, the involution $\iota$ extends to a biregular involution of the blow-up of $Y^{[2]}$ at the image of $\Fano_1(\sfV_5)$ under the map $l\mapsto l\cap Y$. The quotient by this involution is equal to the image $W$ of $\alpha$. It is a EPW-sextic with a unique triple point equal to the image of the exceptional divisor of the blow-up. The EPW-sextic is dual to the EPW-sextic $\calD_Y$.}
\end{itemize}

In our case, the situation is more complicated because $Y$ contains lines and conics. Note that if $\ell_s$ contains a singular point $x$ of $X_{15}$, then a line passing through this point intersects $X_{15}$ at a point on $T(x)$, hence its image is a cycle of two points $(y,\sigma(y))\in E_x\cup \sigma(E_x)$. The surface $\alpha(Y^{[2]})$ must have other singularities.

It is known that the Picard group of $Y^{[2]}$ is generated by $\Pic(Y)$ and $\half [\Delta]$, where $\Delta$ is the exceptional divisor of $\pi:Y^{[2]}\to Y^{(2)}$, the pre-image of the diagonal. Here $\Pic(Y)$ is embedded in $\Pic(Y^{[2]})$ by $D \mapsto \pi^*(p_1^*(D)+p_2^*(D)$, where $p_i:Y\times Y\to Y$ are the projection maps. We denote the image of $D$ in $\Pic(Y^{[2]})$ by $D_2$. Thus we see that $\alpha^*((E_x)_2)$ is the class of a line on $S$ if $x\in \calL$ and the class of a conic if $x\in \calC$. Moreover, $\alpha^*(\half \Delta) = \tilde{B} = -K_S$.

The involution $\sigma$ of $Y_{10}^{[2]}$ extends to an involution $\tilde{\sigma}$ of $Y_{10}^{[2]}$ with fixed locus equal to the image of $S$ in $Y_{10}^{[2]}$ defined by the isomorphism $Y_{10}/(\sigma) = S$. Composing the inclusion $S$ in $Y_{10}^{[2]}$, we get a map
$$f:S\to Y_{10}^{[2]}.$$

\begin{itemize}
\item [(40)]\it{The map $\alpha\circ f:S\to |I_{Y_{10}}(2)|^* = \bbP^5$ coincides with the anti-canonical map of the del Pezzo surface $S$.}  
\end{itemize}

Note that there is also a rational map 
$\Foc(S)\da \Delta = \bbP(\Omega_{Y_{10}}^1) \subset Y_{10}^{[2]}$
that assigns to a nonsingular point of $\Foc(S)$ the cycle $\xi$ that consists of the point and the  infinitely near point corresponding to the  tangent direction defined by the ray passing through this point. We extend this map to $Y_{10}$  by assigning to a point $y\in E_x$ the point $(y,t_y)$, where $t_y$ is the tangent line to the line (conic) $E_x$  at $y$.  Composing this map with the map \eqref{alpha}, we obtain a rational map
$$Y_{10}\to  \bbP(\Omega_{Y_{10}}^1)\da  |I_{2}(Y_{10})|^* \cong \bbP^5.$$

\section{The Picard group}\label{S:picard}
\textbf{8.1} In this section we 
 compute the Picard group of the minimal resolution $Y$ of a \emph{moduli general $15$-nodal quartic surface} $X_{15}$. Since $X_{15}$ varies in a 4-dimensional family, we expect that the rank of $\Pic(Y)$ is equal to 16. Obvious generators over $\bbQ$ are the classes of exceptional curves $E_i$ and the class of a hyperplane section $\eta_H$. 

Recall that a subset $N$ of the set $\calN$ of  nodes of $X_{15}$ is called \emph{even} (resp. \emph{weakly even}) if the sum $e_N$ of the divisor classes of the exceptional curves $E_x,x\in N,$ is divisible by 2 (resp. is divisible by 2 after adding the class $\eta_H$). Let $N$ be the set of $15$ nodes of $X_{15}$ and 
$\calN$ be the sublattice of $\Pic(Y)$ generated by the classes of the exceptional curves $E_x,x\in N,$ and $\eta_H$. We have $V = \calN/2\calN \cong \bbF_2^{16}$ and the kernel of the natural homomorphism $V\to \Pic(Y)/2\Pic(Y)$ is a linear code $C$ generated by the projections of the classes  $e_N, \eta_H$ to $V$. It is known that for a $15$-nodal quartic surface the code $C$ is of dimension $5$ and it contains ten words of weight 6, six words of weight $10$ and $15$ words of weight 8 (see \cite[Theorem 3.3]{Endrass} and \cite[p.48]{Rohn}). After we fix a splitting  $N= \calL+\calC$, formula \eqref{hyperplane} gives an explicit description of $10$  weakly even sets defining words  of weight $6$ and formula \eqref{tropeconic} gives an expression of  weakly even sets defining words of weight $8$.  One chooses a basis of $C$ consisting of $5$ words of weight $6$ from \eqref{hyperplane} by taking 
$B_\calC = (\sigma(E_{12}),\sigma(E_{23}),\sigma(E_{34}),\sigma(E_{45}),\sigma(E_{15})$ (modulo $2\Pic(Y)$).

Let $\tilde{\calN}$ be the even lattice that contains $\calN$ that corresponds to the isotropic  subgroup $H$ of the discriminant group 
$D_{\calN} = \calN^\vee/\calN $ generated by the basis $B_C$ of the code $C$. It contain $\calN$ as a sublattice of finite index and its discriminant group  is isomorphic to  $A = H^\perp/H$ \cite[Proposition 1.4.1]{Nikulin}. It is generated by 
$$\frac{1}{4}\eta_H-\half(E_{14}+E_{25}+E_{35}+E_{56}), \quad  \half(E_{13}+E_{16}+E_{26}+E_{36}),$$
$$\half(E_{13}+E_{25}+E_{34}+E_{56}), \quad \half(E_{13}+E_{24}+E_{12}+E_{46}),$$
$$\half(E_{13}+E_{35}+E_{16}+E_{56}), \quad \half(E_{14}+E_{24}+E_{16}+E_{26}).$$
Comparing this discriminant group and the quadratic form defined on it with the discriminant group of the lattice  $T = U(2)\oplus U(2)\oplus A_1(2)\oplus A_1$, we find that they are isomorphic. It follows from \cite[Theorem 1.4.2]{Nikulin} that the genus of this lattice consists of one isomorphism class.  It also follows from \cite[Theorem 1.4.4]{Nikulin} that the 
lattice $T$ admits a primitive embedding into the K3 lattice $L_{K3} = U^{\oplus 3}\oplus E_8^{\oplus 2}$ and all such embeddings are equivalent with respect to the group of isometries of $L_{K3}$. Let $M$ be the orthogonal complement of $T$ in $L_{K3}$. Then it is an even lattice of signature $(1,15)$ to the discriminant group isomorphic to the discriminant group of $T$ and hence is isomorphic to the discriminant group of $\tilde{\calN}$. 
Again applying Nikulin's result, we find that $\tilde{\calN}\cong M$ and $M$ admits a unique (up to isometry of $L_{K3}$) primitive embedding into $L_{K3}$. The theory of periods of lattice polarized K3 surfaces shows that there exists an irreducible $4$ dimensional moduli space of K3 surfaces containing the lattice $M$ primitively embedded into their Picard lattice. Moreover, the Picard lattice of a moduli general surface isomorphic to $M$.  Since the isomorphism class of $Y$ belongs to this moduli space, a general $15$-nodal quartic surface has the Picard lattice isomorphic to $M$. 

\begin{itemize}
\item[(41)] \it{The transcendental lattice of a K3 surface $Y$ birationally isomorphic to a general 15-nodal quartic surface $X_{15}$ is isomorphic to $U(2)\oplus U(2)\oplus A_1(2)\oplus A_1$. The Picard lattice is isomorphic to the orthogonal complement of this lattice in the K3-lattice $U^{\oplus 3}\oplus  E_8^{\oplus 2}$}.
\item[(42)]\it{The Picard lattice $\Pic(Y)$ is generated by  $15$ classes of the curves $E_x, x\in \calL\cup \calC,$ and 5 classes of the curves 
$\sigma(E_{12}),\sigma(E_{23}),\sigma(E_{34}),$
$\sigma(E_{45},\sigma(E_{15})$.}
\end{itemize} 

It is natural to embed $\Pic(Y)$ into the Picard lattice of a  minimal resolution $Y$ of a Kummer surface $\Kum(J(C))$ of the Jacobian of a general curve $C$ of genus 2. It is known that the latter is generated by the classes of a hyperplane section $\eta_H$, 16 exceptional curves  $N_\alpha$ and 16 classes of trope-conics $T_\beta$ indexed by subsets of cardinality 2 of the set $[1,6]$  with two additional curves $N_0$ and $T_0$. We have $N_\alpha\cdot T_\beta = 1$ if and only if $ \alpha+\beta \in \{(0,(16),(26),(36),(46),(56)\}$, where the addition of subsets is the symmetric sum modulo the complementary set. 

We continue to index our curves $E_x,x\in \calL,$ by 2-element subsets of $[1,5]=\{1,\ldots,5\}$ and $E_x,x\in \calC,$ by elements of $[1,5]$. Then the embedding is as follows:
$$\eta_H\mapsto \eta_H, \quad [E_{ab}] \mapsto [N_{ab}], \quad [E_{a}] \mapsto [N_{a6}],$$
\beq\label{incl}
[\sigma(E_{ab})]\mapsto [T_{ab}], \quad  [\sigma(E_{a})] \mapsto [T_{a6}+T_0+N_0].
\eeq
It is easy to see that under this embedding $\Pic(Y)$ is equal to the orthogonal complement of $r= [N_0]$ with $r^2 = -2$. This embedding corresponds to the fact that, under the specialization, the trope-quartic  $\sigma(E_a)$ splits into the union of two trope-conics $T_{a6}+T_0$ passing through the new node with the exceptional curve $N_0$.

Recall that the group of birational automorphisms  of a general $16$-nodal Kummer quartic surface contains $16$ \emph{switch involutions} which interchange nodes with tropes. One of this involution defines a biregular involution of $Y$ that interchanges $N_\alpha$ with $T_{\alpha}$ \cite{Kondo}. It follows from the above that the involution $\sigma$ is equal to the composition of  this switch involution and the reflection with respect to the $(-2)$-vector 
$\alpha_0 = [T_0+E_0]$.

\textbf{8.2} We know that $X_{15}$ contains $10$ trope-conics and $5$ trope-quartics. They are realized as the images of the exceptional curves $E_x$ under the involution $\sigma$ defined by the choice of a syntheme. Let us see that a general $15$-nodal quartic surface $X_{15}$ has no more conics or quartic rational curves lying on it.

Suppose we have an irreducible  conic $C$ on $X_{15}$. Let $\bar{C}$ be its proper inverse transform on $Y$. Then 
$(\eta_H-\bar{C})^2 = -2, (\eta_H-\bar{C})\cdot \eta_H = 2$, hence $\eta_H = \bar{C}+\bar{C'}$ for some other conic $C'$. Suppose $C = C'$, then the plane $\la C\ra$ is tangent to $X_{15}$ along $C$ and hence $2[\bar{C}] = \eta_H-\sum_{x\in A} E_x$, where $A$ is a set of 6  exceptional curves (to make the self-intersection equal to $-8$). It is easy to see that we can find a trope-conic that passes through 3 nodes in $A$. Computing the intersection number, we conclude that $\bar{C}$ coincides with this conic. Suppose $C\ne C'$. Then $[\bar{C}]+[\bar{C'}]= \eta_H-\sum_{x\in B}E_x$. Since $\Pic(Y)_\bbQ$ is generated by $\eta_H, E_x$, we see that $2[\bar{C}] = H-\sum_{x\in B'}E_x, 2[\bar{C'}] = H-\sum_{x\in B''}E_x$, and we find  a contradiction by computing the intersection with some trope-conic.

Suppose we have an irreducible rational quartic curve $C$ with a node at some singular point $x_0$ on $X_{15}$. Then, as above, we have $\bar{C}\sim \eta_H-E_{x_0}-\half\sum_{i\in I}E_{x_i}$, where $\# I= 8$. Without loss of generality we may assume that $E_{x0} = E_{16}$. Consider 4 trope-conics $\sigma(E_{1a}), a = 2,3,4,5$ passing through $x_0$. Let $Z(a)$ be the set of nodes contained in the conic except the node $x_0$. Let us identify nodes of $X_{15}$ and the exceptional curves on $Y$ with their indices $(ab)$.  We have $Z(a)\cap Z(b) = \{(cd)\}$, where $\{a,b,c,d\} = [1,4]$. It is easy to see that  the set $I$ of 8 nodes has at least 3 nodes in common with one of the subsets $Z(a)$ unless $I = \{(26),(36),(46),(56),(12),(13),(14),(15)\}$. In the latter case $\bar{Q}$ coincides with the quartic curve $\sigma(E_{16})$. So, we have proved the following.

\begin{itemize}
\item[(43)] \it{A general $15$-nodal quartic surface contains $10$ conics which coincide with the trope-conics  $\sigma(E_x),x\in \calL$. It also contains  $15$ rational quartic curves with a double point  at one of the nodes. They coincide with the curves ${}^s\sigma(E_x),x\in \calC$, where ${}^s\sigma$ is one of the six involutions obtained from $\sigma$ by conjugation by a permutation $s\in \frakS_6$.}
\item[(44)]\it{The surface does not contain curves of odd degree.}
\end{itemize}

The last property follows from the fact that $\eta_H$ intersects evenly any divisor class on $Y$ since this is true for generators of $\Pic(Y)$.

\section{Admissible pentads}\label{S:pentads}
\textbf{9.1} 
So far, we have located the following involutions of $X_{15}$:
\begin{itemize}
\item $6$ involutions $\sigma_i$ with quotient a del Pezzo surface of degree 5 with the branch curve a nonsingular curve of genus 6;
\item $6$ Reye involutions $\tau_{\Rey}^i$ corresponding to the six structures of a quartic symmetroid with quotient a Coble surface $Z$ with $K_Z^2 = -5$ and the branch curve equal to the union of five smooth rational curves.
\end{itemize}

Also we  have  obvious involutions defined by the projections from nodes. In this section we will study a new set of involutions defined by  pentads of nodes, no four of the them are coplanar. Such a pentad will be called \emph{admissible}. For a general 15-nodal quartic surface this condition is equivalent to the property that no four of the points lie on a trope-conic.

 Let $\calP =\{x_1,\ldots,x_5\}$ be an admissible pentad. For a general point $x\in X_{15}$, there exists a unique twisted cubic $\gamma$ (maybe degenerate) that passes through $x_1,\ldots,x_5\in \calP$ and $x$. It intersects the surface at  one more point (may be equal to $x$). This defines a birational involution of $X_{15}$ that extends to a biregular involution $\tau_{\calP}$ of $Y$. 
Let $C$ be a quartic curve on $X_{15}$ cut out by a plane   $\Pi$  containing three nodes, say $p_1,p_2,p_3$, from the pentad. Any conic $\gamma'$ in $H$ passing through these nodes intersects $C$ at two additional points $x,x'$. The union of $\gamma'$ and the line $\la p_4,p_5\ra$ is a degenerate twisted cubic $\gamma$ intersecting $X_{15}$ at $p_1,\ldots,p_5$ and $x,x'$. Thus the pair $(x,x')$ is an orbit of the involution $\tau_\calP$ and hence $C$ is invariant curve with respect to this involution.

Note that if we drop the condition of admissibility of a pentad, then a pentad with four nodes on a trope-conic will still define an involution but it will coincide with the involution defined by the projection from the remaining fifth node. If we put all 5 nodes on a trope-conic, the involution is not defined.

In the special case when the plane $H$ cuts out a trope-conic, we see that any  point on it must be fixed under the involution. Also, if there is a trope-quartic $C$ with a double point at $x_1$ and passing through the remaining points of $\calP$, then a twisted cubic $\gamma$ passing through $\calP$ and a general point $x$ on a trope-quartic $C$ will be contained in the cone $K(x_1)$ and hence it will intersect $X_{15}$ at $x$ with multiplicity two. This shows that $C$ is fixed pointwise too.

\begin{itemize} 
\item[(45)] \it{The involution $\tau_{\calP}$ associated with an admissible pentad leaves any plane section $H$ of $X_{15}$ containing three nodes from the pentad invariant. A trope-conic through $3$ points in $\calP$ or a trope-quartic passing through $\calP$ with a double point at one of the points in $\calP$ is fixed pointwise.}
\end{itemize}

Let $\calP = \{x_1,x_2,x_3,x_4,x_5\}$ be an admissible pentad and $\tau_\calP$ be the corresponding involution of $Y$. One checks that the sublattice $M(\calP)$ spanned by the ten divisor classes 
$\eta_H-E_{x_i}-E_{x_j}-E_{x_k}$ is of rank 5. Since $\tau_{\calP}$ obviously leaves invariant the exceptional curves over the remaining 10 nodes, we see that it acts as a reflection on the Picard lattice of $Y$. In fact, we check that the vector 
\beq\label{admisroot}
r_\calP:= 3\eta_H-2\sum_{x\in \calP}E_{x}\in \Pic(Y)
\eeq
has $r_{\calP}^2 = -4$ and it  is orthogonal to the sublattice $M(\calP)$, hence $\tau_\calP^*$ coincides with the reflection with respect to $r_\calP$ and acts as 
\beq\label{admisref}
\tau_\calP^*(v) = v+\half (r_\calP\cdot v)\cdot r_\calP.
\eeq  
In particular, $\tau_{\calP}^*$ leaves invariant the divisor class 
\beq\label{admfiber}
F_i^{\calP} = 2\eta_H-2E_{x_i}-\sum_{j\ne i}E_{x_j}
\eeq 
of self-intersection $0$. Since we can write $F_i^{\calP} = (\eta_H-E_{x_1}-E_{x_2}-E_{x_3})+(\eta_H-E_{x_1}-E_{x_4}-E_{x_5})$, and each of the summand represents a nef class, we see that the divisor class $F_i^{\calP}$ is nef and represents the class of a general fiber of an elliptic pencil $|F_i^{\calP}|$ on $Y$. 

We have $F_i^{\calP}\cdot F_j^{\calP} = 2$ and 
\beq\label{sumfibers}
|\half(F_1^{\calP}+\cdots+F_5^{\calP})| = |5\eta_H-3\sum_{x\in \calP}E_x|.
\eeq

It is clear that ten curves $E_x,x\not\in\calP,$ are irreducible components of fibers of $|F_i^{\calP}|$. 

Suppose   $x_1,x_2,x_3$ and $x_1,x_3,x_4$ do not lie on a trope-conic. Then the plane section of $X_{15}$ containing the points $x_1,x_2,x_3$ or $x_1,x_4,x_5$ cut out two smooth rational curve whose proper transform on $Y$ are linearly equivalent to $\eta_H-E_{x_1}-E_{x_2}-E_{x_3}$ and $\eta_H-E_{x_1}-E_{x_4}-E_{x_5}$. Their sum is linearly  equivalent to $F_1^{\calP}$. This shows that the elliptic pencil $|F_1^{\calP}|$ acquires a reducible fiber of type $\tilde{A}_1$ (or $I_2$ in Kodaira's notation) that consists of two irreducible components (they may be tangent so the fiber becomes of type III in  Kodaira's notation).

On the other hand suppose  that $x_1,x_2,x_3$ lie on a trope-conic passing through $3$ more points $y_1,y_2,y_3\not\in \calP$, and $x_1,x_4,x_5$ do not lie on a trope-conic. Then $F_1^{\calP}$ is linearly equivalent to the sum
{\small $$2(\half(\eta_H-E_{x_1}-E_{x_2}-E_{x_3}-E_{y_1}-E_{y_2}-E_{y_3})+E_{y_1}+E_{y_2}+E_{y_3}+(\eta_H-E_{x_1}-E_{x_4}-E_{x_5}),$$}
where each of the  summands represents the divisor class of a $(-2)$-curve. This shows that the elliptic pencil $|F_1^{\calP}|$  contains a reducible fiber of type $\tilde{D}_4$ (or $I_0^*$ in Kodaira's notation).  

Similar argument shows that if $x_1,x_2,x_3$ and  $x_1,x_4,x_5$  lie on trope-conics, then $|F_1^{\calP}|$ acquires a fiber  of type $\tilde{D}_6$. Also, if $\half(2\eta_H-2E_{x_1}-\sum_{i=2}^5E_{x_i}-\sum_{i=1}^4E_{y_i})$ represents a trope-quartic $C$, then $2C+\sum_{i=1}^4E_{y_i}\in  |F_1^{\calP}|$ is a fiber of type $\tilde{D}_4$.
One can check that no other types of reducible fibers are possible  because $\calP$ is admissible.

Let $\Gamma(\calP)$ be the graph with the set of vertices $[1,6]$ and the set of 5 edges 
corresponding to the indices $(ab)$ of points from $\calP$.  A trope-conic 
 passes through $6$ points corresponding to a duad $(ab)$, three duads in $[1,6]\setminus \{a,b,c\}$  and two duads $(ac),(bc)$. It can be represented by the following graph
\beq\label{graph}
\xy
(0,0)*{\bullet};(10,0)*{\bullet};(20,0)*{\bullet};(40,0)*{\bullet};(50,0)*{\bullet};
(45,-5)*{\bullet};(-10,0)*{\bullet};
(-10, 0)*{};(20,0)*{}**\dir{-};(40, 0)*{};(50,0)*{}**\dir{-};(40, 0)*{};(45,-5)*{}**\dir{-};(50, 0)*{};(45,-5)*{}**\dir{-};
(0,3)*{a};(-10,3)*{c};(10,3)*{b};(20,3)*{c};(40,3)*{d};(50,3)*{e};(48,-5)*{f};
\endxy
\eeq  
By definition, a pentad is admissible if and only if $\Gamma(\calP)$  does not contain a subgraph of this graph with four edges. Also  three nodes in $\calP$ are on a trope-conic if and only if $\Gamma(\calP)$  contain a subgraph of this graph with three  edges.

  The following Table \ref{table:admpentads} shows all possible graphs $\Gamma(\calP)$, the number $c$ of subsets of three nodes on a trope-conic, possible types of reducible fibers of the elliptic pencils $|F_i^{\calP}|$, possible singularities of the branch curve of the double plane realization of $X_{15}$ defined by $\tau_{\calP}$ (see next subsection) and the number of orbits of $\frakS_6$ on the set of pentads (the last two columns were kindly provided by I. Shimada).

\begin{table}
{\footnotesize \begin{equation}
\renewcommand{\arraystretch}{2}
\begin{array}{|c|c|c|c|c|c|}\hline
\textrm{Type}&\textrm{Graph}&c&\textrm{Pencils}&\textrm{Double planes}&$\#$\textrm{orbits}
\\
\hline
1&\raise -12pt \hbox{\includegraphics[scale=.4]{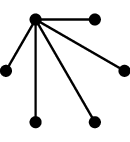}}&0&\parbox{3cm}{$(9\tilde{A}_1+\tilde{D}_4)\times 5$}&\parbox{3cm}{$(14A_1)\times 5$}&6\\
\hline
2&\raise -12pt \hbox{\includegraphics[scale=.4]{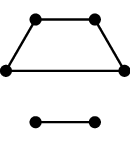}}&4&\parbox{3cm}{$(5\tilde{A}_1+2\tilde{D}_4)\times 4$\\$(\tilde{A}_1+2\tilde{D}_6)\times 1$} &\parbox{3cm}{$(6A_1+2D_4)\times 4$\\$(14A_1)\times 1$} &45\\
\hline
3&\raise -12pt \hbox{\includegraphics[scale=.4]{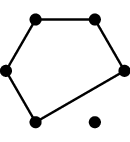}}&5&\parbox{3cm}{$(3\tilde{A}_1+\tilde{D}_4+\tilde{D}_6)\times 5$}&\parbox{3cm}{$(6A_1+2D_4)\times 5$}&72\\
\hline
4&\raise -12pt \hbox{\includegraphics[scale=.4]{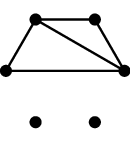}}&2&\parbox{3cm}{$(9\tilde{A}_1+\tilde{D}_4)\times 4$\\$(3\tilde{A}_1+\tilde{D}_4+\tilde{D}_6)\times 1$} &\parbox{3cm}{$(10A_1+D_4)\times 4$\\$(14A_1)\times 1$}&90\\
\hline
5&\raise -12pt \hbox{\includegraphics[scale=.4]{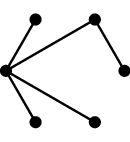}}&3&
\parbox{3cm}{$(5\tilde{A}_1+2\tilde{D}_4)\times 3$\\ $(\tilde{A}_1+3\tilde{D}_4)\times 1$\\ $(9\tilde{A}_1+\tilde{D}_4) \times 1$} &
\parbox{3cm}{$(10A_1+D_4)\times 3$\\ $(14A_1)\times 1$ \\ $(2A_1+3D_4)\times 1$ }&120\\ 
\hline
6&\raise -12pt \hbox{\includegraphics[scale=.4]{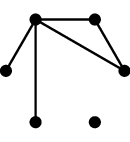}}&2&
\parbox{3cm}{$(5\tilde{A}_1+2\tilde{D}_4)\times 2$\\$(9\tilde{A}_1+\tilde{D}_4)\times 2$\\$(7\tilde{A}_1+\tilde{D}_6)\times 1$}&
\parbox{3cm}{$(10A_1+D_4)\times 4$\\ $(14A_1)\times 1$} &180\\ 
\hline
7&\raise -12pt \hbox{\includegraphics[scale=.4]{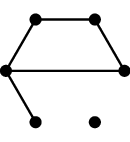}}&3&
\parbox{3cm}{$(5\tilde{A}_1+2\tilde{D}_4)\times 2$\\ $(9\tilde{A}_1+\tilde{D}_4)\times 2$ \\ $(3\tilde{A}_1+\tilde{D}_4+\tilde{D}_6)\times 1$}&
\parbox{3cm}{$(6A_1+2D_4)\times 2$\\ $(10A_1+D_4)\times 2$ \\ $(14A_1)\times 1$}&360\\ 
\hline
8&
\raise -12pt \hbox{\includegraphics[scale=.4]{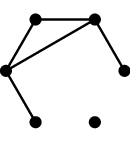}}&3&
\parbox{3cm}{$(5\tilde{A}_1+2\tilde{D}_4)\times 3$\\ $(7\tilde{A}_1+\tilde{D}_6)\times 2$}&
\parbox{3cm}{$(10A_1+D_4)\times 4$\\ $(6A_1+2D_4)\times 1$} &360\\ 
\hline
9&\raise -12pt \hbox{\includegraphics[scale=.4]{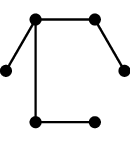}}&4&
\parbox{3cm}{$(5\tilde{A}_1+2\tilde{D}_4)\times 2$\\ $(7\tilde{A}_1+\tilde{D}_6)\times 2$\\$(3\tilde{A}_1+\tilde{D}_6+\tilde{D}_4)\times 1$}&
\parbox{3cm}{$(10A_1+D_4)\times 3$\\
$(6A_1+2D_4)\times 2$}&360\\
\hline
\end{array}
\end{equation}}
\caption{Admissible pentads}\label{table:admpentads}
\end{table}

Note that  an elliptic pencil $|F_i^{\calP}|$ has a section if and only if there exists a trope-conic passing through the node $x_i$ and one other  node of $\calP$. Not all elliptic pencils have a section. In fact one of the fibrations of types 2,4, and 7 has no sections.
Table \ref{table:admpentads} shows that the rank of the sublattice of $\Pic(Y)$ generated by irreducible components of each pencil is equal to $14$. Applying the Shioda-Tate formula, we obtain the following.

\begin{itemize}
\item[(46)] \it{The Mordell-Weyl group of sections of the jacobian fibration of the elliptic fibration defined by $|F_i^{\calP}|$ is of rank 1}.
\end{itemize}

\textbf{9.2}  Consider the quotient $Z$ of $Y$ by the involution $\tau_\calP$. Since any admissible pentad $\calP = \{x_1,\ldots,x_5\}$ consists of $5$ points in general linear position in $\bbP^3$, we can use them to define the map \eqref{kapranov} and identify a twisted cubic through $\calP$ with a point on the irreducible component $F(\sfS_3)_6$ of the surface of lines on the Segre cubic. This defines a $\tau_{\calP}$-equivariant degree $2$ regular map $f:X_{15}\to \dP_5$ to a quintic del Pezzo surface $\dP_5$ that factors through a birational morphism
$$\tilde{f}:Y/(\tau_{\calP}) \to \dP_5.$$
Let $|C_i|, i = 1,\ldots,5,$ be the pencils of conics on $\dP_5$. It is easy to see that in our identification of $\dP_5$ with $F(\sfS_3)_6$, each pencil consists of twisted cubics that are tangent to a ruling of a quadratic cone with vertex at  one of the points $x_i$ and passing through three remaining points $x_j$. If we compare this with our definition of the pencils $|F_i^{\calP}|$ in \eqref{admfiber}, we find that the pre-image of the pencil of conics $|C_i|$ on $\dP_5$ is the pencil $|F_i^{\calP}|$. Since $C_1+\cdots+C_5 \sim -2K_{\dP_5}$, we see that 
the pre-image of a hyperplane section of $\dP_5$ in the anti-canonical embedding is a member of the linear system $|D|$ defined in \eqref{sumfibers}. Comparing the dimensions of the linear systems we find that the linear subsystem $\bbP(H^0(Y,\calO_Y(D))^{\tau_{\calP}})$ of $\tau_{\calP}$-invariant sections is a hyperplane and defines the $\tau_{\calP}$-equivariant map $\tilde{f}$.

\begin{itemize}
\item[(47)] \it{The linear system \eqref{sumfibers} contains a hyperplane that defines a $\tau_{\calP}$-equivariant map  $Y\to \dP_5$. The pre-image of a pencil of conics on $\dP_5$ is an elliptic pencil $|F_i^{\calP}|$ on $Y$.}
\end{itemize} 

Suppose three of the nodes in $\calP$ lie on a trope-conic. Then intersecting the divisor  class of its proper transfer on $Y$ with the divisor class defining the linear system \eqref{sumfibers}, we obtain that its image on the quotient del Pezzo surface is a line in the anti-canonical embedding. Since the trope-conic is contained in the fixed locus of $\tau_{\calP}$,  we see that  the involutions defined by pentads of types $2-9$  are different from the Reye involutions (whose components of the fixed locus are mapped to conics on the del Pezzo surface). Similar computation shows that the image of a trope-quartic  component of the fixed locus is a conic on the  del Pezzo surface.

\begin{itemize}
\item[(48)] \it{The image on $\dP_5$ of a component of the fixed locus of $\tau_{\calP}$ formed by a trope-conic (resp. trope-quartic) is a line (resp. conic) on $\sfD_5$.}
\end{itemize}

The quotient surface $Z = Y/(\tau_{\calP})$ is a smooth rational surface and the map 
$\tilde{f}:Z\to \dP_5$ is a birational morphism. Let $R$ be an exceptional curve that is blown down to a point in $\dP_5$. Then $R$ intersects each $F_i^{\calP}$ with multiplicity zero, hence is contained in a reducible fiber of each elliptic pencil  
$|F_i^{\calP}|$. Of course, the converse is also true. We have $10$ such curves $E_{x},x\not\in \calP$. We can deduce from our description of reducible fibers of the elliptic pencils $|F_i^{\calP}|$ that there are no more such curves, hence $Z$ is obtained from $\dP_5$ by blowing up $10$ points and $K_Z^2 = -5$.

Since the fixed locus $Y^{\tau_{\calP}}$ contains a curve, for example a multiple component of a fiber of type $\tilde{D}_4$, the involution  acts as the minus identity on the transcendental part of $H^2(Y,\bbC)$. This shows that the invariant part of the involution $\tau_\calP^*$ on $H^2(Y,\bbC)$ coincides with the invariant part on $\Pic(Y)$. Since $\tau_{\calP}$ acts as a reflection on $\Pic(Y)$, the dimension of the invariant part $\Pic(Y)_\bbC^{\tau_{\calP}}$ is equal to $15$, and therefore the Lefschetz number $\textrm{Lef}(\tau_{\calP})$ is equal to $15-7+2 = 10$.  Applying the Lefschetz fixed-point-formula, we obtain that the Euler-Poincar\'e characteristic of the fixed locus $Y^{\tau_\calP}$ of $\tau_\calP$ is equal to $10$. It follows that the fixed locus  contains at least $5$ disjoint smooth rational curves $R_1,\ldots,R_5$. Using the Hurwitz formula, we also see that the Euler-Poincar\'e characteristic of $Z = Y/(\tau_{\calP})$ is equal to $17$. This confirms our previous assertion that $Z$ is obtained from $\dP_5$ by blowing up $10$ points.

\begin{itemize}
\item[(49)] \it{The quotient surface $Z = Y/(\tau_\calP)$ is a rational  surface $Z$ with $K_Z^2 = -5$. The branch curve $B$ of the cover $Y\to Z$ belongs to $|-2K_Z|$ and contains $5$ disjoint smooth rational curves with self-intersection $-4$. The surface $Z$  is isomorphic to the blow-up of 10 points on a quintic del Pezzo surface, the points are the images of nodes on $X_{15}$ which do not belong to $\calP$. The image of $B$ on $Z$ belongs to $|-2K_{\dP_5}|$}. 
\end{itemize}

 We will see in Example \ref{non-coble} that the fixed locus may contain other non-rational component and the quotient surface  $Z$ may not be a Coble surface as in the case of a Reye involution.

\begin{remark}Recall that a quintic del Pezzo surface $\dP_5$ admits $5$ non-projectively equivalent blow-down morphism $\pi:\dP_5\to \bbP^2$. They are defined by the linear system $|-K_{\dP_5}-C_i|$, where $|C_i|$ is pencil of conics. Since we know that the pre-image of $|-K_{\dP_5}|$ on $Y$ is  linear system \eqref{sumfibers} and the pre-image of $|C_i|$ is the elliptic pencil $|F_i^{\calP}|$, we obtain that the linear system
\beq\label{doubleplane}
|3\eta_H-E_{x_i}-2\sum_{j\ne i}E_{x_j}|
\eeq
defines a morphism $Y\to \dP_5\to \bbP^2$. Thus we have five double plane realizations of $Y$ corresponding to a choice of an admissible pentad. The branch curve of these double planes are reducible plane sextics and the possible types of their singularities  are given in Table \ref{table:admpentads}.

We have to warn that the double plane with the branch curve described in  Table \ref{table:admpentads} is not isomorphic in general to a $15$-nodal quartic surface. In  fact, formula \eqref{doubleplane} shows that the  set of the images of the curves $E_{x_i},x_i\in \calP,$ on the del Pezzo quotient surface are smooth rational curves $R_i$ of degree $6$ not passing through the singular points of the branch curve. We have $R\cdot K_{\dP_5} = -6$, hence $R_i^2 = 4$. In order that $R_i$ be the image of a $(-2)$-curve on $Y$, it must split under the cover. So we need to find $6$ such curves which are everywhere tangent to the branch curve $B$. It is easy to see that a smooth rational sextic of $\dP_5$ is equal to the proper transform under the blow-up $\dP_5\to \bbP^2$ of either a smooth conic not passing through the fundamental points or a quartic curve with nodes at three of the fundamental points. In both cases the class of $R_i$ is divisible by $2$. As soon as we have found such a set of curves $R_1,\ldots,R_5$ everywhere tangent to the branch curve,  we use formula \eqref{admfiber} to define the linear system $|\eta_H|$ that will map the double cover to a $X_{15}$.
\end{remark}

\begin{example} Let $\calP$ be an admissible pentad of type 1. This is the only type such that all subsets of 3 nodes in $\calP$ do not lie on a trope-conic. We call such pentad a \emph{G\"opel pentad} (by analogy with a G\"opel tetrad on a Kummer surface \cite[Chapter VII, \S 53]{Hudson}). We may assume that it coincides with the set $\calC$ which we index by the set $\{(16),(26),(36),(46),(56)\}$. The fixed locus of $\tau_\calP$ is equal to the union of the trope-quartics  $\sigma(E_x),x\in \calC$. The surface $Y$ is the double cover of the blow-up of 10 intersection points of the intersection points of $5$ conics on a quintic del Pezzo surface. The  branch curve is equal to the proper inverse transform of these conics on the blow-up.  We recognize a similar construction of $X_{15}$ from a Reye involution $\tau_{\Rey}$. In that case, the fixed locus was the union of the curves $E_y,y\in \calC$. The six involutions defined by G\"opel pentads  are in fact conjugate to the six Reye involutions
$$\tau_{\calP_i} = \sigma_i\circ \tau_{\Rey_i}\circ\sigma_i, \quad i = 1,\ldots,6.$$ 
\end{example}

\begin{example}  Assume $\calP$ is of type $2$. We can index its points by 
$(12),(13),$ and $(24),(34),(56)$. There are 4 trope-conics passing through  three nodes from $\calP$. They are $\sigma(E_{15}), \sigma(E_{25}), \sigma(E_{35}), \sigma(E_{45})$. The node $(56)$ is contained in all trope-conics. Each pair of them 
passing through two complementary pairs of the remaining nodes form a fiber of type $\tilde{D}_6$ of $F_{56}$. The remaining fibers have two reducible fibers of type $\tilde{D}_4$. 

Let us consider the images of the four trope-conics on $\dP_5$. We observe that $E_{14}$ and $E_{23}$ are mapped to the intersection points of the two disjoint pairs of lines $\{\ell_1,\ell_2|$ and $\{\ell_3,\ell_4\}$, the images of $\sigma(E_{25}), \sigma(E_{35}$ and the images of $\sigma(E_{15}),\sigma(E_{45}$. The remaining $8$ exceptional curves are blown down to $8$ points, two on each of the lines. Since the image of the branch curve $B$ on $Z$ belongs to $|-2K_{\dP_5}|$, we see that $B$ has another component $C$ of degree $6$ in the anti-canonical embedding. Let us consider a  $\dP_5$ as the blow-up of points $p_1,p_2,p_3,p_4$ in $\bbP^2$. It is easy to see that the images of the four lines $\ell_1,\ldots,\ell_4$ are either two pairs of lines 
joining two complementary pairs of points or two lines $\la p_i,p_j\ra$ and $\la p_i,p_k\ra$ intersecting at one of the points $p_i$. There will be $4$ non-projectively equivalent  ways to blow down to obtain the first possibility and only one to  obtain the second possibility. In the first case the remaining irreducible component $C$ of $B$ is the proper transform of a nonsingular conic not passing through $p_1,\ldots,p_4$. In the second case, the curve $C$ is the proper transform  of a quartic curve with three nodes at $p_j,p_k$ and $p_l$, where $\{i,j,k,l\} = \{1,2,3,4\}$. Thus we see that $|-2K_Z|$ consists of $5$ disjoint smooth rational curves and the surface $Z$ is a Coble surface.
\end{example}

\begin{example}\label{non-coble} Let us consider  a  pentad of type 3. By (48), the image of the fixed locus of $\tau_{\rey}$ on $\dP_5$ contains five lines among its irreducible components. 
They are represented by a pentagon subgraph of the Petersen graph. There are five projective equivalence classes of blowing down morphisms $Z\to \dP_5$. The image of the  pentagon of lines is formed by three lines $\la p_i,p_j\ra, \la p_i,p_k\ra, \la p_i,p_l\ra$. The remaining component of the image of the fixed locus is the proper transform of a cubic curve through the four points. It is a quintic elliptic curve on $\dP_5$. Thus we see that the fixed locus of $\tau_{\calP}$ may contain an elliptic curve besides  five $(-2)$-curves. Since the pentagon of lines belongs to $|-K_{\dP_5}|$ and passes through all ten points which we blow-up to obtain $Z$, we see that $|-K_Z|\ne \emptyset$ and \emph{$Z$ it is not a Coble surface}.
\end{example}

 \end{document}